\documentclass[twoside]{article}

\usepackage[accepted]{aistats2023}
\usepackage[round]{natbib}

\usepackage{mathtools}
\usepackage{appendix}
\usepackage{amsthm}
\usepackage{amsmath}

\usepackage{amsfonts}
\usepackage{multirow}
\usepackage{multicol}

\usepackage{color}

\usepackage{xcolor,colortbl}
\definecolor{LightBlue}{rgb}{0.78,1,1}

\usepackage{amsfonts}
\usepackage{multirow}
\usepackage{tcolorbox}

\usepackage{algorithm}
\usepackage{algorithmic}

\usepackage{graphicx}
\usepackage{subfigure}
\usepackage{amsmath,amssymb,amsthm,wrapfig,enumitem}

\newtheorem{theorem}{Theorem}
\newtheorem{lemma}{Lemma}

\newtheorem{assumption}{Assumption}
\newtheorem{remark}{Remark}
\newtheorem{corollary}{Corollary}

\begin{document}

%
\runningtitle{AdaGDA: Faster Adaptive Gradient Descent Ascent Methods for Minimax Optimization}

%
\runningauthor{Feihu Huang, Xidong Wu, Zhengmian Hu}

\twocolumn[

\aistatstitle{ AdaGDA: Faster Adaptive Gradient Descent Ascent Methods \\
for Minimax Optimization }

\aistatsauthor{ Feihu Huang$^{1,2,*}$ \And Xidong Wu$^3$ \And  Zhengmian Hu$^3$ }

\aistatsaddress{ 1. College of Computer Science and Technology, Nanjing University of Aeronautics and Astronautics, Nanjing, China; \\ 2. MIIT Key Laboratory of Pattern Analysis and Machine Intelligence, China; *E-mail: huangfeihu2018@gmail.com;
\\  3. Department of Electrical and Computer Engineering, University of Pittsburgh, Pittsburgh, USA. } ]

\begin{abstract}
  In the paper, we propose a class of faster adaptive Gradient Descent Ascent (GDA) methods for solving the nonconvex-strongly-concave minimax problems by using the unified adaptive matrices,
 which include almost all existing coordinate-wise and global adaptive learning rates.
 In particular, we provide an effective convergence analysis framework for our adaptive GDA methods.
  Specifically, we propose a fast Adaptive Gradient Descent Ascent (AdaGDA) method based on the basic momentum technique, which reaches a lower gradient complexity of $\tilde{O}(\kappa^4\epsilon^{-4})$ for finding an $\epsilon$-stationary point without large batches, which improves the existing results of the adaptive GDA methods by a factor of $O(\sqrt{\kappa})$. Moreover, we propose an accelerated version of AdaGDA (VR-AdaGDA) method based on the momentum-based variance reduced technique, which achieves a lower gradient complexity of $\tilde{O}(\kappa^{4.5}\epsilon^{-3})$ for finding an $\epsilon$-stationary point without large batches, which improves the existing results of the adaptive GDA methods by a factor of $O(\epsilon^{-1})$. Moreover, we prove that
  our VR-AdaGDA method can reach the best known gradient complexity of $\tilde{O}(\kappa^{3}\epsilon^{-3})$ with the mini-batch size $O(\kappa^3)$. The experiments on policy evaluation and fair classifier learning tasks are conducted to verify the efficiency of our new algorithms.
\end{abstract}

\section{Introduction}

In the paper, we consider the following stochastic nonconvex-strongly-concave minimax problem:
\begin{align} \label{eq:1}
 \min_{x \in \mathcal{X} } \max_{y \in \mathcal{Y}} \ \mathbb{E}_{\xi \sim \mathcal{D}}[f(x,y;\xi)],
\end{align}
where function $f(x,y)=\mathbb{E}_{\xi}[f(x,y;\xi)]: \mathbb{R}^{d_1}\times\mathbb{R}^{d_2} \rightarrow \mathbb{R}$ is $\mu$-strongly concave over $y$ but possibly nonconvex over $x$,
and $\xi$ is a random variable following an unknown distribution $\mathcal{D}$.
Here $\mathcal{X}\subseteq \mathbb{R}^{d_1}$ and $ \mathcal{Y}\subseteq \mathbb{R}^{d_2}$ are  nonempty compact convex sets.
In fact, Problem \eqref{eq:1} is widely used to many machine learning applications, such as adversarial training \citep{goodfellow2014generative,tramer2018ensemble,nouiehed2019solving}, reinforcement learning \citep{wai2019variance}
and robust federated learning \citep{deng2021distributionally}. In the following, we specifically provide two popular applications
that can be formulated as the above Problem \eqref{eq:1}.

\begin{table*}
  \centering
  \caption{ \textbf{Gradient complexity} comparison of the representative gradient descent ascent methods for finding
  an $\epsilon$-stationary point of the \textbf{nonconvex-strongly-concave} problem \eqref{eq:1}, i.e., $\mathbb{E}\|\nabla F(x)\|\leq \epsilon$
  or its equivalent variants, where $F(x)=\max_{y\in \mathcal{Y}}f(x,y)$. \textbf{ALR} is adaptive learning rate.  \textbf{Cons$(x,y)$} denotes constraint sets on variables x and y, respectively.
  Here \textbf{Y} denotes the fact that there exists a convex constraint set on variable, otherwise is \textbf{N}.
  \textbf{1} denotes Lipschitz continuous of $\nabla_xf(x,y)$, $\nabla_yf(x,y)$ for all $x,y$; \textbf{2} means Lipschitz continuous of $\nabla_xf(x,y;\xi)$, $\nabla_yf(x,y;\xi)$ for all $\xi, x,y$; \textbf{3} denotes the bounded set $\mathcal{Y}$ with a diameter $D\geq 0$. Since some algorithms do not provide the explicit dependence on $\kappa$, we use $p(\kappa)$. }
  \label{tab:1}
   \resizebox{\textwidth}{!}{
\begin{tabular}{c|c|c|c|c|c|c|c}
  \hline
   \textbf{Algorithm} & \textbf{Reference}  & \textbf{Cons$(x,y)$} & \textbf{Loop(s)} & \textbf{Batch Size} & \textbf{Complexity}& \textbf{ALR} & \textbf{Conditions}  \\ \hline
  SGDA  & \cite{lin2019gradient} & N,\ Y & Single & $O(\kappa\epsilon^{-2})$ & $O(\kappa^3\epsilon^{-4})$ & & \textbf{1}, \textbf{3} \\  \hline
  SREDA  & \cite{luo2020stochastic} & N,\ Y & Double& $O(\kappa^2\epsilon^{-2})$ & $O(\kappa^3\epsilon^{-3})$ & & \textbf{2} \\  \hline
  Acc-MDA  & \cite{huang2022accelerated} & Y,\ Y & Single & $O(1)$ & $\tilde{O}(\kappa^{4.5}\epsilon^{-3})$ & & \textbf{2} \\  \hline
  Acc-MDA  & \cite{huang2022accelerated} & Y,\ Y & Single & $O(\kappa^3)$ & $\tilde{O}(\kappa^{3}\epsilon^{-3})$ & & \textbf{2} \\  \hline
  PDAda & \cite{guo2021novel} & N,\ Y & Single & $O(1)$ & $O(\kappa^{4.5}\epsilon^{-4})$ & $\surd$ & \textbf{1} \\  \hline
  NeAda-AdaGrad & \cite{yang2022nest} & N,\ Y & Double & $O(\epsilon^{-2})$ & $\tilde{O}(p(\kappa)\epsilon^{-4})$ & $\surd$ & \textbf{1} \\  \hline
 AdaGDA & Ours & Y,\ Y & Single & {\color{red}{ $O(1)$ }} & {\color{red}{ $\tilde{O}(\kappa^{4}\epsilon^{-4})$ }} & $\surd$ & \textbf{1} \\ \hline
 VR-AdaGDA & Ours & Y,\ Y & Single & {\color{red}{ $O(1)$ }} & {\color{red}{ $\tilde{O}(\kappa^{4.5}\epsilon^{-3})$}} & $\surd$  & \textbf{2}  \\ \hline
 VR-AdaGDA & Ours & Y,\ Y & Single & {\color{red}{ $O(\kappa^3)$}} & {\color{red}{ $\tilde{O}(\kappa^{3}\epsilon^{-3})$ }} & $\surd$ & \textbf{2} \\ \hline
\end{tabular}
 }
\end{table*}

\textbf{1) Policy Evaluation.} Policy evaluation
aims at estimating the value function corresponding to a certain policy, which is a stepping stone of policy
optimization and serves as an essential component of many reinforcement learning algorithms such as actor-critic algorithm \citep{konda2000actor}.
Specifically, we consider a Markov decision process (MDP) $(\mathcal{S}, \mathcal{A}, \mathcal{P}, R, \tau)$,
where $\mathcal{S}$ denotes the state space, and $\mathcal{A}$ denotes the action
space, and $\mathcal{P}(s'|s,a)$ denotes the transition kernel to the
next state $s'$ given the current state $s$ and action $a$,
and $\tau \in [0,1)$ is the discount factor. $R(s,a,s') \in [-r,r] \ (r>0)$ is an immediate reward once an
agent takes action $a$ at state $s$ and transits to state $s'$, and $R(s,a)$ is the reward at $(s,a)$, defined as
$R(s,a) = \mathbb{E}_{s'\sim \mathcal{P}(\cdot|s,a)}[R(s,a,s')]$.
$\pi(s,a) : \mathcal{S}\times \mathcal{A}\rightarrow \mathbb{R}$ denotes a stationary
policy that is the probability of taking action $a\in \mathcal{A}$ given the current state $s\in \mathcal{S}$.
We let $V^{\pi}(s) = \mathbb{E}\big[ \sum_{t=0}^{+\infty}\tau^t R(s_t,a_t)|s_0=s,\pi \big]$ denote state value function.
Further let $V(s;\theta)$ be the parameterized approximate function of $V^{\pi}(s)$,
and $V(s;\theta)$ generally is a smooth nonlinear function.
Following \cite{wai2019variance}, we can solve the following minimax problem to find an optimal approximated value function,
defined as
\begin{align} \label{eq:2}
 \min_{\theta \in \Theta} \max_{\omega \in \mathbb{R}^d} & \ \mathbb{E}_{s,a,s'}\bigg[  \langle\delta \nabla_{\theta}V(s;\theta),\omega \rangle  \nonumber \\
& \quad - \frac{1}{2}\omega^T \big(\nabla_{\theta}V(s;\theta) \nabla_{\theta}V(s;\theta)^T\big) \omega \bigg],
\end{align}
where $\delta = R(s,a,s') + \tau V_{\theta}(s') - V_{\theta}(s)$, and $\mathbb{E}_{s,a,s'}$ is taking expectation for $s\sim d^{\pi}(\cdot)$ that is stationary distribution of states, $a\in \pi(\cdot,s)$ and $s'\sim \mathcal{P}(\cdot|s,a)$. Here matrix $H_{\theta}= \mathbb{E}\big[\nabla_{\theta}V(s;\theta) \nabla_{\theta}V(s;\theta)^T\big]$ is generally positive definite. The above problem \eqref{eq:2} is generally nonconvex on variable $\theta$ when using the neural networks to approximate value function $V^{\pi}(s)$.

\textbf{2) Robust Federated Averaging.}
Federated Learning (FL)~\citep{mcmahan2017communication} is a popular learning paradigm for training a centralized model
based on decentralized data over a network of clients.
Specifically, we have $n$ clients in FL framework, and $\mathcal{D}_i$ is the data distribution on $i$-th device,
and the data distributions $\{\mathcal{D}_i\}_{i=1}^n$ generally are different. The goal of FL is to
learn a global variable $w$ based on these heterogeneous data from different data distributions.
To well solve the data heterogeneity issue in FL,
some robust FL methods \citep{deng2021distributionally,reisizadeh2020robust} have been proposed,
which solve the following distributionally robust empirical loss problem:
\begin{align}
 \min_{w \in \Omega} \max_{p\in \Pi} \bigg\{\sum_{i=1}^n p_i \mathbb{E}_{\xi \sim \mathcal{D}_i}[f_i(w;\xi)] - \lambda\psi(p) \bigg\},
\end{align}
where $p_i\in (0,1)$ denotes the proportion of $i$-th device in the entire model, and $f_i(w;\xi)$ is
the loss function on $i$-th device, and $\lambda>0$ is a tuning parameter, and $\psi(p)$
is a (strongly) convex regularization. Here $\Pi = \{ p\in \mathbb{R}^n : \sum_{i=1}^np_i=1, \ p_i\geq 0\}$ is
a $n$-dimensional simplex, and $\Omega \subseteq \mathbb{R}^d$ is a nonempty convex set.

Since the above minimax problem \eqref{eq:1} frequently appeared in many machine learning applications,
multiple methods have been proposed to solve it. For example,  \cite{lin2019gradient,lin2020near}
proposed a stochastic gradient descent ascent (SGDA) method to solve the problem \eqref{eq:1}.
Subsequently, a class of accelerated SGDA methods \citep{luo2020stochastic,huang2022accelerated} have been presented  based on the variance reduced techniques of
SPIDER \citep{fang2018spider,wang2019spiderboost} and STORM \citep{cutkosky2019momentum}, respectively.
More recently, \cite{guo2021novel,yang2022nest} introduced the adaptive versions of SGDA by using the adaptive learning rates. However, these adaptive SGDA methods still suffer from the high sample (gradient) complexities (please see Table~\ref{tab:1}). Meanwhile, the adaptive PDAda algorithm in \cite{guo2021novel} only considers using adaptive learning rate in updating minimized variable $x$.
Thus, there exists a natural question:
\begin{center}
\begin{tcolorbox}
\textbf{ Can we develop faster adaptive gradient descent ascent methods to solve the Problem \eqref{eq:1}, which use adaptive learning rates in updating
both variables $x$ and $y$ ? }
\end{tcolorbox}
\end{center}

In the paper, we give an affirmative answer to the above question and propose
a class of faster adaptive gradient descent ascent methods to solve the Problem \eqref{eq:1}.
Our methods can use many types of adaptive learning rates to update both variables $x$ and $y$.
Moreover, our methods can flexibly incorporate momentum and variance-reduced techniques.
Our main contributions are in three-fold:
\begin{itemize}
\item[(1)] We propose a class of faster adaptive gradient descent ascent methods for the nonconvex-strongly-concave minimax Problem \eqref{eq:1} using the universal adaptive matrices for both variables $x$ and
$y$, which include most existing adaptive learning rates.
\item[(2)] We propose a fast adaptive gradient descent ascent (AdaGDA) method based on the basic momentum technique used in Adam algorithm \citep{kingma2014adam}.
Meanwhile, we present an accelerated version of AdaGDA (VR-AdaGDA) method based on the momentum-based variance reduced technique used in STORM algorithm \citep{cutkosky2019momentum}.
\item[(3)] We provide an effective convergence analysis framework for our adaptive methods under mild assumptions. Specifically, we prove that our AdaGDA method has a gradient complexity of $\tilde{O}(\kappa^4\epsilon^{-4})$ without large batches, which improves the existing result of adaptive method for solving the problem \eqref{eq:1} by a factor of $O(\kappa^{1/2})$. Our VR-AdaGDA method has a lower gradient complexity of $\tilde{O}(\kappa^{4.5}\epsilon^{-3})$ without large batches, which improves the existing best known result by a factor of $O(\epsilon^{-1})$ (please see Table~\ref{tab:1} for comparison summary).
\end{itemize}

From Table~\ref{tab:1}, despite achieving a better rate when compared to PDAda~\citep{guo2021novel} and
NeAda-AdaGrad~\citep{yang2022nest}, our VR-AdaGDA algorithm still have the same complexity rate as the existing non-adaptive Acc-MDA algorithm.
In fact, only under some specific cases such as
sparse gradient condition, the adaptive gradient methods have a
faster convergence rate than the non-adaptive counterparts.
For example, Adagrad~\citep{duchi2011adaptive} shows a better convergence rate than SGD under the sparse gradient condition.
In fact, we propose an adaptive gradient-based algorithm framework for minimax optimization based on the general adaptive matrices without
some specific conditions such as sparse gradients. It is well known that adaptive gradient methods generally perform
well in practice although with same convergence rate as non-adaptive gradient methods.
In fact, our VR-AdaGDA algorithm obtains a near-optimal complexity $O(\epsilon^{-3})$
in finding an $\epsilon$-stationary point (i.e., $\mathbb{E}||\nabla F(x)||\leq \epsilon$, where $F(x)=\max_{y}\mathbb{E}[f(x,y;\xi)]$). Thus, we can not obtain a lower complexity than this near-optimal complexity $\tilde{O}(\epsilon^{-3})$.
\textbf{NOTE THAT:} the single-level problem
\begin{align}
\min_{x\in \mathbb{R}^d} f(x)\equiv\mathbb{E}_{\xi}[f(x;\xi)] \label{eq:3}
\end{align}
can be seen as a specific case of the minimax Problem (\ref{eq:1}).
For example, $f(x,y;\xi)=af(x;\xi)+b$, where $a>0$ and $b\geq 0$ are constants, i.e., given any $x$,
the function $f(x,\cdot;\xi)=c$ is independent on $x$ and $\xi$, where $c$ is a constant.
\cite{arjevani2019lower} proves the stochastic algorithms in solving the single-level nonconvex stochastic problem (\ref{eq:3})
has a lower bound complexity $O(\epsilon^{-3})$ for finding an $\epsilon$-stationary point (i.e., $\mathbb{E}||\nabla f(x)||\leq \epsilon$).
Since the above Problem (\ref{eq:3}) can be seen as a specific case of the minimax Problem (\ref{eq:1}),
the stochastic algorithms in solving the minimax stochastic Problem (\ref{eq:1})
also has a lower bound complexity $O(\epsilon^{-3})$ for finding an $\epsilon$-stationary point (i.e., $\mathbb{E}||\nabla F(x)||\leq \epsilon$).

\section{Related Works}
In this section, we overview the existing first-order methods for minimax optimization and
adaptive gradient methods.

\subsection{ Minimax Optimization Methods }
Minimax optimization has recently been shown great successes in many machine learning applications
such as adversarial training, robust federated learning, and policy optimization.
Thus, many first-order methods \citep{nouiehed2019solving,lin2019gradient,lin2020near,lu2020hybrid,yan2020optimal,yang2020catalyst,yang2020global,rafique2021weakly,liu2021first} were recently proposed to solve the minimax problems.
For example, some (stochastic) gradient-based descent ascent methods \citep{lin2019gradient,nouiehed2019solving,lu2020hybrid,yan2020optimal,lin2020near}
have been proposed for solving the minimax problems.
Subsequently, several accelerated gradient descent ascent algorithms \citep{rafique2021weakly,luo2020stochastic,huang2022accelerated} were proposed to solve the stochastic minimax problems based on the variance-reduced techniques.
Meanwhile, \cite{huang2021efficient,chen2021proximal} studied the nonsmooth nonconvex-strongly-concave minimax optimization.
In addition, \cite{huang2022accelerated,wang2022zeroth} studied the zeroth-order methods for solving the nonconvex-strongly-concave minimax problems.
\cite{huang2023gradient} have proposed a class of Riemanian gradient descent ascent algorithms to
 solve the geodesically-nonconvex strongly-concave minimax problems on
 Riemanian manifolds.
\cite{zhang2021complexity,li2021complexity} studied the lower bound complexities of nonconvex-strongly-concave minimax optimization.
More recently, \cite{guo2021novel,yang2022nest} proposed an adaptive gradient descent ascent method for solving Problem \eqref{eq:1}.

\subsection{ Adaptive Gradient Methods }
Adaptive gradient methods are a class of popular optimization tools to solve large-scale machine learning problems, \emph{e.g.}, Adam \citep{kingma2014adam}
is one of the most popular optimization tools for training deep neural networks (DNNs),
which is a version of the first adaptive gradient method, AdaGrad~\citep{duchi2011adaptive}.
The adaptive gradient methods have been widely studied in machine learning community. Among them, Adam \citep{kingma2014adam} is the most popular one and uses a coordinate-wise adaptive learning rate
and momentum technique to accelerate algorithm.
Multiple variants of Adam algorithm \citep{reddi2019convergence,chen2018convergence,guo2021novel}
have been presented to obtain a convergence guarantee under the nonconvex setting.
Due to the coordinate-wise adaptive learning rate, Adam often shows a bad generalization performance in training DNNs.
To improve the generalization performance of Adam, recently several adaptive gradient methods such as AdamW \citep{loshchilov2017decoupled} and AdaBelief \citep{zhuang2020adabelief} were developed.
More recently, the accelerated adaptive gradient methods \citep{cutkosky2019momentum,huang2021super} were designed based on the variance-reduced techniques.
In particular, \cite{huang2021super} proposed a faster and universal adaptive gradient SUPER-ADAM framework using a universal adaptive matrix.

\subsection{Notations}
For vectors $x$ and $y$, $x^r \ (r>0)$ denotes the element-wise
power operation, $x/y$ denotes the element-wise division and $\max(x,y)$ denotes the element-wise maximum. $I_{d}$ denotes a $d$-dimensional identity matrix. For two vectors $x$ and $y$, $\langle x,y\rangle$ is their inner product. $\|\cdot\|$ denotes
the $\ell_2$ norm for vectors and spectral norm for matrices, respectively.
$\nabla_x f(x,y)$ and $\nabla_y f(x,y)$
are the partial derivatives \emph{w.r.t.} variables $x$ and $y$ respectively.  $I_d$ denotes $d$-dimension identity matrix.
$a=O(b)$ means that $a\leq C b$ for some constant $C>0$, and the notation $\tilde{O}(\cdot)$ hides logarithmic terms.
Given the mini-batch samples $\mathcal{B}=\{\xi_i\}_{i=1}^q$, we let $\nabla f(x;\mathcal{B})=\frac{1}{q}\sum_{i=1}^q \nabla f(x;\xi_i)$.

\section{ Faster Adaptive Gradient Descent Ascent Methods }
In this section, we propose a class of faster adaptive gradient descent ascent methods for
solving the minimax problem \eqref{eq:1}.
Specifically, we propose a fast adaptive gradient descent ascent (AdaGDA) based on the basic momentum technique of Adam \citep{kingma2014adam}.
Meanwhile, we further propose an accelerated version of AdaGDA (VR-AdaGDA) based on the momentum-based variance reduced technique of STORM \citep{cutkosky2019momentum}.

\subsection{ AdaGDA Algorithm }
We first propose a new fast adaptive gradient descent ascent (AdaGDA) algorithm for
solving the Problem \eqref{eq:1} based on the basic momentum technique.
Algorithm \ref{alg:1} summarizes the algorithmic framework of our AdaGDA.

\begin{algorithm}[tb]
\caption{ AdaGDA Algorithm }
\label{alg:1}
\begin{algorithmic}[1] 
\STATE {\bfseries Input:} $T$, tuning parameters $\{\gamma, \lambda, \eta_t, \alpha_t, \beta_t\}_{t=1}^T$ and mini-batch size $q$; \\
\STATE {\bfseries initialize:} Initial input $x_1 \in \mathcal{X}$, $y_1 \in \mathcal{Y}$, and draw a mini-batch i.i.d. samples $\mathcal{B}_1=\{\xi_i^1\}_{i=1}^q$,
and then compute $v_1=\nabla_x f(x_1,y_1;\mathcal{B}_1)=\frac{1}{q}\sum_{i=1}^q\nabla_x f(x_1,y_1;\xi^1_i)$ and $w_1 = \nabla_y f(x_1,y_1;\mathcal{B}_1)=\frac{1}{q}\sum_{i=1}^q\nabla_y f(x_1,y_1;\xi^1_i)$; \\
\FOR{$t = 1, 2, \ldots, T-1$}
\STATE Generate the adaptive matrices $A_t \in \mathbb{R}^{d_1 \times d_1}$ and $B_t \in \mathbb{R}^{d_2 \times d_2}$; \\
\textcolor{blue}{One example: $A_t$ and $B_t$ are generated from (\ref{eq:4}) and (\ref{eq:5}), respectively.} \\
\STATE $x_{t+1} = x_t + \eta_t (\tilde{x}_{t+1}-x_t)$ with $\tilde{x}_{t+1} = \arg\min_{x \in \mathcal{X}} \big\{\langle v_t, x\rangle + \frac{1}{2\gamma}(x-x_t)^TA_t(x-x_t) \big\} $;  \\
\STATE $y_{t+1} = y_t + \eta_t (\tilde{y}_{t+1}-y_t) $ with $\tilde{y}_{t+1} = \arg\max_{y \in \mathcal{Y}} \big\{\langle w_t, y\rangle - \frac{1}{2\lambda}(y-y_t)^TB_t(y-y_t)\big\} $; \\
\STATE Draw a mini-batch i.i.d. samples $\mathcal{B}_{t+1}=\{\xi^{t+1}_i\}_{i=1}^q$, and then compute
\STATE $v_{t+1} = \alpha_{t+1}\nabla_x f(x_{t+1},y_{t+1};\mathcal{B}_{t+1}) + (1-\alpha_{t+1})v_t $;
\STATE $w_{t+1} = \beta_{t+1}\nabla_y f(x_{t+1},y_{t+1};\mathcal{B}_{t+1}) + (1-\beta_{t+1})w_t $;
\ENDFOR
\STATE {\bfseries Output:} $x_{\zeta}$ and $y_{\zeta}$ chosen uniformly random from $\{x_t, y_t\}_{t=1}^{T}$.
\end{algorithmic}
\end{algorithm}

At the line 4 of Algorithm \ref{alg:1}, we generate the adaptive matrices $A_t$ and $B_t$
for variables $x$ and $y$, respectively.
Specifically, we use the general adaptive matrix $A_t\succeq \rho I_{d_1}$ for variable $x$ as in the SUPER-ADAM \citep{huang2021super},
and the global adaptive matrix $B_t = b_t I_{d_2} \ (b_t>0)$.
For example, we can generate the matrix $A_t$ as in the Adam \citep{kingma2014adam}, defined as:
\begin{align}
 & \tilde{v}_0=0, \ \tilde{v}_t = \varrho \tilde{v}_{t-1} + (1-\varrho) \nabla_x f(x_t,y_t;\xi_t)^2, \nonumber \\
 & A_t=\mbox{diag}( \sqrt{\tilde{v}_t} + \rho), \ t\geq 1,  \label{eq:4}
\end{align}
where $\varrho \in (0,1)$ and $\rho >0$.
Matrix $B_t$ is defined as: given $\beta \in (0,1)$ and $\varrho >0$,
\begin{align}
  & b_0>0, \ b_t = \varrho b_{t-1} + (1-\varrho)\|\nabla_y f(x_t,y_t;\xi_t)\|, \nonumber \\
  & B_t= (b_t + \rho) I_{d_2}, \  t\geq 1, \label{eq:5}
\end{align}
 which can be seen as a new global adaptive learning rate.
 Meanwhile, we also generate the matrix $A_t$ as in the AdaBelief \citep{zhuang2020adabelief}, defined as:
\begin{align}
 & \tilde{v}_0=0, \ \tilde{v}_t = \varrho \tilde{v}_{t-1} + (1-\varrho) \big(\nabla_x f(x_t,y_t;\xi_t)-v_t\big)^2, \nonumber \\
 & A_t=\mbox{diag}( \sqrt{\tilde{v}_t} + \rho), \ t\geq 1,
\end{align}
where $\varrho \in (0,1)$ and $\rho >0$.
Matrix $B_t$ is defined as:
\begin{align}
  & b_0>0, \ b_t = \varrho b_{t-1} + (1-\varrho)\|\nabla_y f(x_t,y_t;\xi_t)-w_t\|, \nonumber \\
  & B_t= (b_t + \rho) I_{d_2}, \  t\geq 1,
\end{align}
where $\varrho \in (0,1)$ and $\rho >0$.

At the lines 5 and 6 of Algorithm \ref{alg:1}, we apply the generalized projection gradient iteration to update
variables $x$ and $y$ based on the adaptive matrices $A_t$ and $B_t$, respectively. Meanwhile, we use the momentum iteration to update the variables $x$ and $y$.
At the lines 8 and 9 of Algorithm \ref{alg:1}, we adopt the basic momentum technique to estimate the
stochastic gradients $v_t$ and $w_t$.

\subsection{ VR-AdaGDA Algorithm}
Next, we propose an accelerated version of AdaGDA (VR-AdaGDA) algorithm
based on the momentum-based variance reduced technique.
Algorithm \ref{alg:2} shows the algorithmic framework of the VR-AdaGDA.

\begin{algorithm}[tb]
\caption{ VR-AdaGDA Algorithm }
\label{alg:2}
\begin{algorithmic}[1] 
\STATE {\bfseries Input:} $T$, tuning parameters $\{\gamma, \lambda, \eta_t, \alpha_t, \beta_t\}_{t=1}^T$ and mini-batch size $q$;  \\
\STATE {\bfseries initialize:} Initial input $x_1 \in \mathcal{X}$, $y_1 \in \mathcal{Y}$, and draw a  mini-batch i.i.d. samples $\mathcal{B}_1=\{\xi_i^1\}_{i=1}^q$,
and then compute $v_1=\nabla_x f(x_1,y_1;\mathcal{B}_1)$ and $w_1 = \nabla_y f(x_1,y_1;\mathcal{B}_1)$; \\
\FOR{$t = 1, 2, \ldots, T-1$}
\STATE Generate the adaptive matrices $A_t \in \mathbb{R}^{d_1 \times d_1}$ and $B_t \in \mathbb{R}^{d_2 \times d_2}$; \\
\textcolor{blue}{One example: $A_t$ and $B_t$ are generated from (\ref{eq:4}) and (\ref{eq:5}), respectively.}
\STATE $x_{t+1} = x_t + \eta_t (\tilde{x}_{t+1}-x_t)$ with $\tilde{x}_{t+1} = \arg\min_{x \in \mathcal{X}} \big\{\langle v_t, x\rangle + \frac{1}{2\gamma}(x-x_t)^TA_t(x-x_t) \big\} $;  \\
\STATE $y_{t+1} = y_t + \eta_t (\tilde{y}_{t+1}-y_t) $ with $\tilde{y}_{t+1} = \arg\max_{y \in \mathcal{Y}} \big\{\langle w_t, y\rangle - \frac{1}{2\lambda}(y-y_t)^TB_t(y-y_t)\big\} $; \\
\STATE Draw a mini-batch i.i.d.  samples $\mathcal{B}_{t+1}=\{\xi^{t+1}_i\}_{i=1}^q$, and then compute
\STATE $v_{t+1} = \nabla_x f(x_{t+1},y_{t+1};\mathcal{B}_{t+1}) + (1-\alpha_{t+1})\big(v_t - \nabla_x f(x_t,y_t;\mathcal{B}_{t+1})\big)  $;
\STATE $w_{t+1} = \nabla_y f(x_{t+1},y_{t+1};\mathcal{B}_{t+1}) + (1-\beta_{t+1})\big(w_t - \nabla_y f(x_t,y_t;\mathcal{B}_{t+1})\big)$;
\ENDFOR
\STATE {\bfseries Output:} $x_{\zeta}$ and $y_{\zeta}$ chosen uniformly random from $\{x_t, y_t\}_{t=1}^{T}$.
\end{algorithmic}
\end{algorithm}

At the lines 5 and 6 of Algorithm \ref{alg:2}, we simultaneously use the momentum iteration and the  generalized projection gradient iteration to update
variables $x$ and $y$.
At the lines 8 and 9 of Algorithm \ref{alg:2}, we apply the momentum-based variance reduced technique to
estimate the stochastic gradients $v_t$ and $w_t$.
For example, the estimator of gradient $\nabla f(x_{t+1},y_{t+1})$ is defined as:
\begin{align}
 v_{t+1}
& = \alpha_{t+1}\nabla_x  f(x_{t+1},y_{t+1};\mathcal{B}_{t+1}) + (1-\alpha_{t+1}) \big[v_t \nonumber \\
& \quad + \nabla_x  f(x_{t+1},y_{t+1};\mathcal{B}_{t+1}) \!-\! \nabla_x  f(x_t,y_t;\mathcal{B}_{t+1})\big]. \nonumber
\end{align}
Compared with the estimator $v_{t+1}$ in Algorithm \ref{alg:1}, $v_{t+1}$ in Algorithm \ref{alg:2}
adds the term $(1-\alpha_{t+1})\big(\nabla f(x_{t+1},y_{t+1};\mathcal{B}_{t+1}) - \nabla f(x_t,y_t;\mathcal{B}_{t+1})\big)$ to reduce variance of gradient estimator, where $\alpha_{t+1}\in (0,1)$.

\section{Convergence Analysis}
In this section, we study the convergence properties of our new algorithms (\emph{i.e.}, AdaGDA and VR-AdaGDA) under mild assumptions.
All related proofs are provided in the following Appendix.

\subsection{Mild Assumptions }
We have the following mild assumptions for Problem \eqref{eq:1}.

\begin{assumption} \label{ass:1}
Each component function $f(x,y;\xi)$ has an unbiased stochastic gradient with
bounded variance $\sigma^2$, i.e., for all $\xi, x \in \mathcal{X}, y \in \mathcal{Y}$, $\mathbb{E}[\nabla_x f(x,y;\xi)] = \nabla_x f(x,y)$,
$\mathbb{E}\|\nabla_x f(x,y) - \nabla_x f(x,y;\xi)\|^2 \leq \sigma^2$, $\mathbb{E}[\nabla_y f(x,y;\xi)] = \nabla_y f(x,y)$
and $\mathbb{E}\|\nabla_y f(x,y) - \nabla_y f(x,y;\xi)\|^2 \leq \sigma^2$.
\end{assumption}

\begin{assumption} \label{ass:2}
Function $f(x,y)$ is $\mu$-strongly concave in $y\in \mathcal{Y}$, i.e., for all $x\in \mathcal{X}$ and $y_1, y_2\in \mathcal{Y}$,
we have $\|\nabla_yf(x,y_1)-\nabla_yf(x,y_2)\| \geq \mu \|y_1-y_2\|$.
Then the following inequality holds
\begin{align}
 f(x,y_1) \leq f(x,y_2) \!+\! \langle\nabla_y f(x,y_2), y_1-y_2\rangle \!-\! \frac{\mu}{2}\|y_1-y_2\|^2. \nonumber
\end{align}

\end{assumption}
Since the function $f(x,y)$ is strongly concave in $y\in \mathcal{Y}$, there exists a unique solution to
the problem $\max_{y\in \mathcal{Y}} f(x,y)$ for any $x$. Here we let $y^*(x) = \arg\max_{y \in \mathcal{Y}} f(x,y)$
and $F(x) = f(x,y^*(x)) = \max_{y\in \mathcal{Y}} f(x,y)$.

\begin{assumption} \label{ass:3}
The function $F(x)$ is bounded below in $\mathcal{X}$, \emph{i.e.,} $F^* = \inf_{x\in \mathcal{X}}F(x) > -\infty$.
\end{assumption}

\begin{assumption} \label{ass:4}
In our algorithms, the adaptive matrices $A_t$ for all $t\geq 1$ for updating the variables $x$ satisfies
$A_t^T = A_t$ and $ \lambda_{\min}(A_t) \geq \rho >0$, where $\rho$ is an  appropriate positive number.
\end{assumption}
Assumption \ref{ass:4} ensures that the adaptive matrices $A_t$ for all $t\geq 1$ are positive definite as in \cite{huang2021super}.
Since the function $f(x,y)$ is $\mu$-strongly concave in $y$, we can easily obtain the global solution of the subproblem $\max_{y\in \mathcal{Y}}f(x,y)$. Without loss of generalization, in the following convergence analysis,
we consider the adaptive matrices $B_t=b_tI_{d_2}$ for all $t\geq 1$ for updating the variables $y$ satisfies
$\hat{b} \geq b_t \geq b>0$, as the global adapitve learning rates \citep{li2019convergence,ward2019adagrad,huang2021super}.

\begin{assumption} \label{ass:5}
The objective function $f(x,y)$ has a $L_f$-Lipschitz gradient, i.e.,
 for all $x,x_1,x_2\in \mathcal{X}$ and $y,y_1,y_2 \in \mathcal{Y}$, we have
\begin{align}
& \|\nabla_x f(x_1,y)-\nabla_x f(x_2,y)\| \leq L_f \|x_1 - x_2\|, \nonumber \\
& \|\nabla_x f(x,y_1)-\nabla_x f(x,y_2)\| \leq L_f \|y_1 - y_2\|, \nonumber \\
& \|\nabla_y f(x_1,y)-\nabla_y f(x_2,y)\| \leq L_f \|x_1 - x_2\|, \nonumber \\
& \|\nabla_y f(x,y_1)-\nabla_y f(x,y_2)\| \leq L_f \|y_1 - y_2\|.  \nonumber
\end{align}
\end{assumption}

\subsection{ Convergence Metrics }
We introduce useful convergence metrics to measure convergence of our algorithms.
Let $\phi_t(x)=\frac{1}{2}x^TA_t x$, according to Assumption 4,
$\phi_t(x)$ is $\rho$-strongly convex. We define a prox-function (\emph{i.e.}, Bregman distance) associated with $\phi_t(x)$ as in \cite{censor1981iterative,censor1992proximal,ghadimi2016mini}:
\begin{align}
 D_t(x,x_t) & = \phi_t(x) - \big[ \phi_t(x_t) + \langle\nabla \phi_t(x_t), x-x_t\rangle\big] \nonumber \\
 & = \frac{1}{2}(x-x_t)^TA_t(x-x_t).
\end{align}
The line 5 of Algorithms \ref{alg:1} or \ref{alg:2} is equivalent to the following generalized projection problem:
\begin{align}
 \tilde{x}_{t+1} = \arg\min_{x\in \mathcal{X}}\big\{ \langle v_t, x\rangle + \frac{1}{\gamma}D_t(x,x_t)\big\}.
\end{align}
As in \cite{ghadimi2016mini}, we define a generalized projected gradient $\mathcal{G}_{\mathcal{X}}(x_t,v_t,\gamma)=\frac{1}{\gamma}(x_t-\tilde{x}_{t+1})$.
At the same time, we define a gradient mapping $\mathcal{G}_{\mathcal{X}}(x_t,\nabla F(x_t),\gamma)=\frac{1}{\gamma}(x_t-x^*_{t+1})$, where
\begin{align}
 x^*_{t+1} = \arg\min_{x\in \mathcal{X}}\big\{ \langle \nabla F(x_t), x\rangle + \frac{1}{\gamma}D_t(x,x_t)\big\}.
\end{align}
For Problem \eqref{eq:1}, when $\mathcal{X}\subset \mathbb{R}^{d_1}$, we use the standard gradient mapping metric $\mathbb{E}\|\mathcal{G}_{\mathcal{X}}(x_t,\nabla F(x_t),\gamma)\|$ to  measure the convergence of
our algorithms, as in \cite{ghadimi2016mini}. When $\mathcal{X} = \mathbb{R}^{d_1}$, we use the standard gradient metric
$\mathbb{E}\|\nabla  F(x_t)\|$ to measure convergence of
our algorithms, as in \cite{lin2019gradient}.

\subsection{Convergence Analysis of the AdaGDA Algorithm}
We analyze the convergence properties of our AdaGDA algorithm under Assumptions \ref{ass:1}, \ref{ass:2}, \ref{ass:3}, \ref{ass:4} and  \ref{ass:5}. The following theorems show our main theoretical results.
The detail proofs are provided in the Appendix \ref{Appendix:A1}.
For notational simplicity, let $L=L_f(1+\kappa)$ and $\kappa=\frac{L_f}{\mu}$.

\begin{theorem} \label{th:1}
Suppose the sequence $\{x_t,y_t\}_{t=1}^T$ be generated from Algorithm \ref{alg:1}. When $\mathcal{X} \subset \mathbb{R}^{d_1}$, and given $B_t=b_tI_{d_2} \ (\hat{b} \geq b_t \geq b>0)$ for all $t\geq 1$,  $\eta_t=\frac{k}{(m+t)^{1/2}}$ for all $t\geq 0$, $\alpha_{t+1}=c_1\eta_t$, $\beta_{t+1}=c_2\eta_t$, $m\geq \max\big(k^2, (c_1k)^2, (c_2k)^2\big)$, $k>0$, $\frac{9\mu^2}{4} \leq c_1 \leq \frac{m^{1/2}}{k}$, $\frac{75L^2_f}{2} \leq c_2 \leq \frac{m^{1/2}}{k}$, $0< \gamma \leq \min\big( \frac{15\sqrt{2}\lambda\mu^2\rho}{2\sqrt{400L^2_f\lambda^2+24\mu^2\lambda^2+16875\hat{b}^2\kappa^2L^2_f\mu^2}}, \frac{m^{1/2}\rho}{4Lk}\big)$ and $0< \lambda \leq \min\big(\frac{405bL^2_f\mu^{3/2}}{8\sqrt{50L^2_f+9\mu^2}},\frac{b}{6L_f}\big)$, we have
\begin{align}
 & \frac{1}{T} \sum_{t=1}^T\mathbb{E}\|\mathcal{G}_{\mathcal{X}}(x_t,\nabla F(x_t),\gamma)\| \nonumber \\
 & \leq \frac{2\sqrt{3G}m^{1/4}}{T^{1/2}} + \frac{2\sqrt{3G}}{T^{1/4}},
\end{align}
where $G = \frac{F(x_1) - F^*}{k\gamma\rho} + \frac{9b_1L^2_f\Delta^2_1}{k\lambda\mu\rho^2} + \frac{2\sigma^2}{q k\mu^2\rho^2} + \frac{2m\sigma^2}{q k\mu^2\rho^2}\ln(m+T)$ and $\Delta^2_1=\|y_1-y^*(x_1)\|^2$.
\end{theorem}

\begin{theorem} \label{th:2}
Assume that the sequence $\{x_t,y_t\}_{t=1}^T$ be generated from the Algorithm \ref{alg:1}. When $\mathcal{X} = \mathbb{R}^{d_1}$, and given $B_t=b_tI_{d_2} \ (\hat{b} \geq b_t \geq b>0)$ for all $t\geq 1$,  $\eta_t=\frac{k}{(m+t)^{1/2}}$ for all $t\geq 0$, $\alpha_{t+1}=c_1\eta_t$, $\beta_{t+1}=c_2\eta_t$, $m\geq \max\big(k^2, (c_1k)^2, (c_2k)^2\big)$, $k>0$, $\frac{9\mu^2}{4} \leq c_1 \leq \frac{m^{1/2}}{k}$, $\frac{75L^2_f}{2} \leq c_2 \leq \frac{m^{1/2}}{k}$, $0< \gamma \leq \min\big(\frac{15\sqrt{2}\lambda\mu^2\rho}{2\sqrt{400L^2_f\lambda^2+24\mu^2\lambda^2+16875\hat{b}^2\kappa^2L^2_f\mu^2}}, \frac{m^{1/2}\rho}{4Lk}\big)$ and $0<\lambda\leq \min\big(\frac{405bL^2_f\mu^{3/2}}{8\sqrt{50L^2_f+9\mu^2}},\frac{b}{6L_f}\big)$, we have
\begin{align}
 & \frac{1}{T} \sum_{t=1}^T \mathbb{E}\|\nabla  F(x_t)\| \nonumber \\
 & \leq \frac{\sqrt{\frac{1}{T}\sum_{t=1}^T\mathbb{E}\|A_t\|^2}}{\rho}\bigg(\frac{2\sqrt{3G'}m^{1/4}}{T^{1/2}} + \frac{2\sqrt{3G'}}{T^{1/4}}\bigg),
\end{align}
where $G' =\frac{\rho(F(x_1) - F^*)}{k\gamma} + \frac{9b_1L^2_f\Delta^2_1}{k\lambda\mu} + \frac{2\sigma^2}{q k\mu^2} + \frac{2m\sigma^2}{q k\mu^2}\ln(m+T)$.
\end{theorem}

\begin{remark}
Without loss of generality, let $k=O(1)$, $b=O(1)$, $\hat{b}=O(1)$ and $\frac{15\sqrt{2}\lambda\mu^2\rho}{2\sqrt{400L^2_f\lambda^2+24\mu^2\lambda^2+9375\hat{b}^2\kappa^2L^2_f\mu^2}} \leq \frac{m^{1/2}\rho}{4Lk}$,
we have $m\geq \max\big(k^2, (c_1k)^2, (c_2k)^2, \frac{225L^2k^2\lambda^2\mu^4}{800L^2_f\lambda^2 +48\mu^2\lambda^2 + 18750\hat{b}^2\kappa^2L^2_f\mu^2} \big)$. At the same time, let $\frac{b}{6L_f} \leq \frac{405bL^2_f\mu^{3/2}}{8\sqrt{50L^2_f+9\mu^2}}$, we have $0< \lambda \leq \frac{b}{6L_f}$. Given $\gamma = \frac{15\sqrt{2}\lambda\mu^2\rho}{2\sqrt{400L^2_f\lambda^2+24\mu^2\lambda^2+9375\hat{b}^2\kappa^2L^2_f\mu^2}}$,
$\lambda =\frac{b}{6L_f}$, $c_1=\frac{9\mu^2}{4}$ and $c_2=\frac{75L^2_f}{2}$.
 Without loss of generality, let $\mu\leq \frac{1}{L_f}$,
 it is easily verified that $\gamma = O(\frac{1}{\kappa^2})$, $\lambda=O(\frac{1}{L_f})$, $c_1=O(\mu^2)$,
$c_2 = O(L^2_f)$. Then we have $m=O(L^4_f)$.  When mini-batch size $q=O(1)$, we have $G=O(\kappa^2 + \kappa^2 \ln(m+T))=\tilde{O}(\kappa^2)$.
Thus, our AdaGDA algorithm has a convergence rate of $\tilde{O}(\frac{\kappa}{T^{1/4}})$. Let $\tilde{O}(\frac{\kappa}{T^{1/4}}) \leq \epsilon$, \emph{i.e.}, $\mathbb{E}\|\mathcal{G}_{\mathcal{X}}(x_\zeta,\nabla F(x_\zeta),\gamma)\|\leq \epsilon$ or $\mathbb{E}\|\nabla  F(x_\zeta)\|\leq \epsilon$, we have $T\leq \kappa^4\epsilon^{-4}$. In Algorithm \ref{alg:1}, we need to compute $2q$ stochastic gradients to estimate partial derivative estimators $v_t$ and $w_t$ at each iteration,
and need $T$ iterations.
Therefore, our AdaGDA algorithm has a gradient (\emph{i.e.}, stochastic first-order oracle) complexity of
$2q\cdot T=\tilde{O}(\kappa^4 \epsilon^{-4})$ for finding an $\epsilon$-stationary point. Note that the term $\sqrt{\frac{1}{T}\sum_{t=1}^T\mathbb{E}\|A_t\|^2}$ is
bounded to the existing adaptive learning rates in Adam algorithm \citep{kingma2014adam} and so on.
For example, given the above adaptive learning rate \eqref{eq:4} and the standard bounded gradient $\|\nabla_x f(x,y)\|\leq \delta$ as in
Adam, we have $\sqrt{\frac{1}{T}\sum_{t=1}^T\mathbb{E}\|A_t\|^2} \leq \delta+\sigma+\rho_0$.
\end{remark}

\subsection{Convergence Analysis of the VR-AdaGDA Algorithm}
We further study the convergence properties of our VR-AdaGDA algorithm under Assumptions \ref{ass:1}, \ref{ass:2}, \ref{ass:3}, \ref{ass:4} and \ref{ass:6}.
The detail proofs are provided in the Appendix \ref{Appendix:A2}.
Here we first use the following assumption instead of the above Assumption \ref{ass:5}.

\begin{assumption} \label{ass:6}
Each component function $f(x,y;\xi)$ has a $L_f$-Lipschitz gradient, i.e.,
 for all $x,x_1,x_2\in \mathcal{X}$ and $y, y_1,y_2 \in \mathcal{Y}$, we have
\begin{align}
& \|\nabla_x f(x_1,y;\xi) - \nabla_x f(x_2,y;\xi)\| \leq L_f \|x_1 - x_2\|, \nonumber \\
& \|\nabla_x f(x,y_1;\xi) - \nabla_x f(x,y_2;\xi)\| \leq L_f \|y_1 - y_2\|, \nonumber \\
& \|\nabla_y f(x_1,y;\xi) - \nabla_y f(x_2,y;\xi)\| \leq L_f \|x_1 - x_2\|, \nonumber \\
& \|\nabla_y f(x,y_1;\xi) - \nabla_y f(x,y_2;\xi)\| \leq L_f \|y_1 - y_2\|. \nonumber
\end{align}
\end{assumption}
By using convexity of $\|\cdot\|$ and Assumption \ref{ass:6}, we have $\|\nabla_x f(x_1,y) - \nabla_x f(x_2,y)\|=\|\mathbb{E}\big[\nabla_x f(x_1,y;\xi) - \nabla_x f(x_2,y;\xi)\big]\| \leq \mathbb{E} \|\nabla_x f(x_1,y;\xi) - \nabla_x f(x_2,y;\xi)\|\leq L_f\|x_1-x_2\|$. Similarly, we also have $\|\nabla_x f(x,y_1) - \nabla_y f(x,y_1)\|\leq L_f\|y_1-y_2\|$, $\|\nabla_y f(x,y_1) - \nabla_y f(x,y_1)\|\leq L_f\|y_1-y_2\|$ and $\|\nabla_y f(x_1,y) - \nabla_y f(x_2,y)\|\leq L_f\|x_1-x_2\|$. In the other words, Assumption \ref{ass:6} includes Assumption \ref{ass:5}, \emph{i.e.}, Assumption \ref{ass:6} is stricter than Assumption \ref{ass:5}.

\begin{theorem} \label{th:3}
Suppose the sequence $\{x_t,y_t\}_{t=1}^T$ be generated from Algorithm \ref{alg:2}. When $\mathcal{X} \subset \mathbb{R}^{d_1}$, and given $B_t=b_tI_{d_2} \ (\hat{b} \geq b_t \geq b>0)$ for all $t\geq 1$,  $\eta_t=\frac{k}{(m+t)^{1/3}}$ for all $t\geq 0$, $\alpha_{t+1}=c_1\eta^2_t$, $\beta_{t+1}=c_2\eta^2_t$, $c_1 \geq \frac{2}{3k^3} + \frac{9\mu^2}{4}$ and $c_2 \geq \frac{2}{3k^3} + \frac{75L^2_f}{2}$, $m \geq \max\big( k^3, (c_1k)^3, (c_2k)^3\big)$,
$0< \lambda \leq \min\big(\frac{27\mu bq}{32},\frac{b}{6L_f}\big)$ and
$0< \gamma \leq \min\big( \frac{\rho\lambda\mu\sqrt{q}}{L_f\sqrt{32\lambda^2+150q\kappa^2\hat{b}^2}}, \frac{m^{1/3}\rho}{2Lk} \big)$,
we have
\begin{align}
 & \frac{1}{T}\sum_{t=1}^T\mathbb{E}\|\mathcal{G}_{\mathcal{X}}(x_t,\nabla F(x_t),\gamma)\| \nonumber \\
 & \leq \frac{2\sqrt{3M}m^{1/6}}{T^{1/2}} + \frac{2\sqrt{3M}}{T^{1/3}},
\end{align}
where $M = \frac{F(x_1) - F^*}{T\gamma k\rho} + \frac{9L^2_f b_1}{k\lambda\mu\rho^2}\Delta_1^2  + \frac{2\sigma^2m^{1/3}}{k^2q\mu^2\rho^2} + \frac{2k^2(c_1^2+c_2^2)\sigma^2}{q\mu^2\rho^2}\ln(m+T)$.
\end{theorem}

\begin{theorem} \label{th:4}
Suppose the sequence $\{x_t,y_t\}_{t=1}^T$ be generated from Algorithm \ref{alg:2}. When $\mathcal{X}=\mathbb{R}^{d_1}$, and  given $B_t=b_tI_{d_2} \ (\hat{b} \geq b_t \geq b>0)$ $\eta_t=\frac{k}{(m+t)^{1/3}}$, $\alpha_{t+1}=c_1\eta^2_t$, $\beta_{t+1}=c_2\eta^2_t$, $c_1 \geq \frac{2}{3k^3} + \frac{9\mu^2}{4}$ and $c_2 \geq \frac{2}{3k^3} + \frac{75L^2_f}{2}$, $m \geq \max\big( k^3, (c_1k)^3, (c_2k)^3\big)$, $0< \lambda \leq \min\big(\frac{27\mu bq}{32},\frac{b}{6L_f}\big)$ and
$0< \gamma \leq \min\big( \frac{\rho\lambda\mu\sqrt{q}}{L_f\sqrt{32\lambda^2+150q\kappa^2\hat{b}^2}}, \frac{m^{1/3}\rho}{2Lk} \big)$,
we have
\begin{align}
 & \frac{1}{T} \sum_{t=1}^T \mathbb{E}\|\nabla F(x_t)\| \nonumber \\
 & \leq \frac{\sqrt{\frac{1}{T}\sum_{t=1}^T\mathbb{E}\|A_t\|^2}}{\rho}\bigg(\frac{2\sqrt{3M'}m^{1/6}}{T^{1/2}} + \frac{2\sqrt{3M'}}{T^{1/3}}\bigg),
\end{align}
where $M' =\frac{\rho(F(x_1) - F^*)}{T\gamma k} + \frac{9L^2_f b_1}{k\lambda\mu}\Delta_1^2  + \frac{2\sigma^2m^{1/3}}{k^2q\mu^2} + \frac{2k^2(c_1^2+c_2^2)\sigma^2}{q\mu^2}\ln(m+T)$.
\end{theorem}

\begin{remark}
Without loss of generality, let $k=O(1)$, $b=O(1)$, $\hat{b}=O(1)$ and $\frac{\rho\lambda\mu\sqrt{q}}{L_f\sqrt{32\lambda^2+150q\kappa^2\hat{b}^2}}\leq \frac{m^{1/3}\rho}{2Lk}$,
we have $m\geq \big(k^3, (c_1k)^3, (c_2k)^3, \frac{8(Lk\lambda\mu)^3q^{3/2}}{L_f(32\lambda^2+150q\kappa^2\hat{b}^2)^{3/2}} \big)$. Given
$\gamma = \frac{\rho\lambda\mu\sqrt{q}}{L_f\sqrt{32\lambda^2+150q\kappa^2\hat{b}^2}}=\frac{\rho\lambda\sqrt{q}}{\kappa\sqrt{32\lambda^2+150q\kappa^2\hat{b}^2}}$
and $\lambda = \min\big(\frac{27\mu bq}{32},\frac{b}{6L_f}\big)$. Without loss of generality, let $\mu\leq \frac{1}{L_f}$, we have $\lambda=O(b\mu)$.
When mini-batch size $q=O(1)$, it is easy to verify that $\gamma=O(\kappa^{-3})$, $\lambda=O(\mu)$, $c_1=O(\mu^2)$, $c_2=O(L^2_f)$ and $m=O(L^6_f)$.
Then we have $M=O\big(\kappa^3 + \kappa +  \kappa^2 +  \kappa^2\ln(m+T)\big) = O(\kappa^3)$. Thus, our VR-AdaGDA algorithm has a convergence rate of $O(\frac{\kappa^{3/2}}{T^{1/3}})$.
Let $O(\frac{\kappa^{3/2}}{T^{1/3}}) \leq \epsilon$,
i.e., $\mathbb{E}\|\mathcal{G}_{\mathcal{X}}(x_\zeta,\nabla F(x_\zeta),\gamma)\|\leq \epsilon$
or $\mathbb{E}\|\nabla  F(x_\zeta)\|\leq \epsilon$, we have $T\leq \kappa^{4.5}\epsilon^{-3}$. In Algorithm \ref{alg:2},
we need to compute $4q$ stochastic gradients to estimate the partial derivative estimators $v_t$ and $w_t$ at each iteration, and need $T$ iterations.
Therefore, our VR-AdaGDA algorithm has a gradient complexity of $4q\cdot T=O(\kappa^{4.5} \epsilon^{-3})$ for finding an $\epsilon$-stationary point.
\end{remark}

\begin{corollary} \label{cor:1}
Under the same conditions of Theorem \ref{th:2}, given mini-batch size $q=O(\kappa^{\nu})$ for $\nu>0$ and $\frac{27\mu bq}{32} \leq \frac{b}{6L_f}$,
i.e., $q=\kappa^{\nu} \leq \frac{16}{81L_f\mu}$, our VR-AdaGDA algorithm has a lower gradient complexity of
$\tilde{O}\big(\kappa^{(4.5-\frac{\nu}{2})}\epsilon^{-3}\big)$ for finding an $\epsilon$-stationary point.
\end{corollary}

\begin{remark}
Without loss of generality, let $\nu=1$, we have $q=\kappa = \frac{L_f}{\mu} \leq \frac{16}{81L_f\mu}$.
Thus, we have $L_f \leq \frac{4}{9}$. Although the objective function $f(x,y)$ in the minimax problem \eqref{eq:1}
maybe not satisfy this condition $L_f\leq \frac{4}{9}$,
we can easily change the original objective function $f(x,y)$ to a new function $\tilde{f}(x,y)=\beta f(x,y), \ \beta>0$. Since $\nabla \tilde{f}(x,y) = \beta \nabla f(x,y)$,
the gradient of function $\tilde{f}(x,y)$ is $\hat{L}$-Lipschitz continuous ($\hat{L}=\beta L_f$).
Thus, we can choose a suitable parameter $\beta$ to ensure this new objective function $\tilde{f}(x,y)$ satisfies the condition $\hat{L}=\beta L_f\leq \frac{4}{9}$.
\end{remark}

\section{Experimental Results}
In this section, we show the empirical results to validate the efficiency of our algorithms on two tasks: 1) Policy Evaluation, and 2) Fair Classifier.
We compare our algorithms (AdaGDA and VR-AdaGDA) with the existing state-of-the-art algorithms in Table \ref{tab:1} for solving nonconvex-strongly-concave minimax problems.

The experiments are run on CPU machines with 2.3 GHz Intel Core i9 as well as NVIDIA Tesla P40 GPU.

\begin{figure*}[ht]
\centering
\subfigure[CartPole-v1]{
\hspace{0pt}
\includegraphics[width=.31\textwidth]{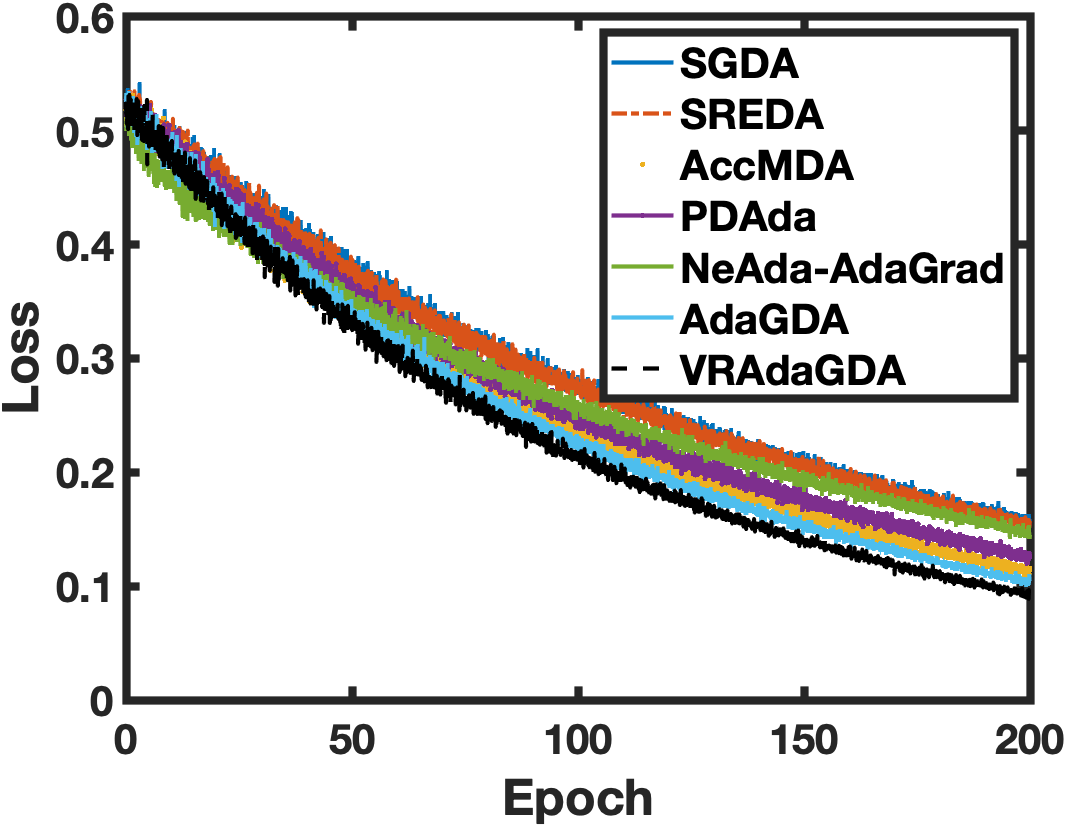}
}
\subfigure[ Acrobat-v1]{
\hspace{0pt}
\includegraphics[width=.31\textwidth]{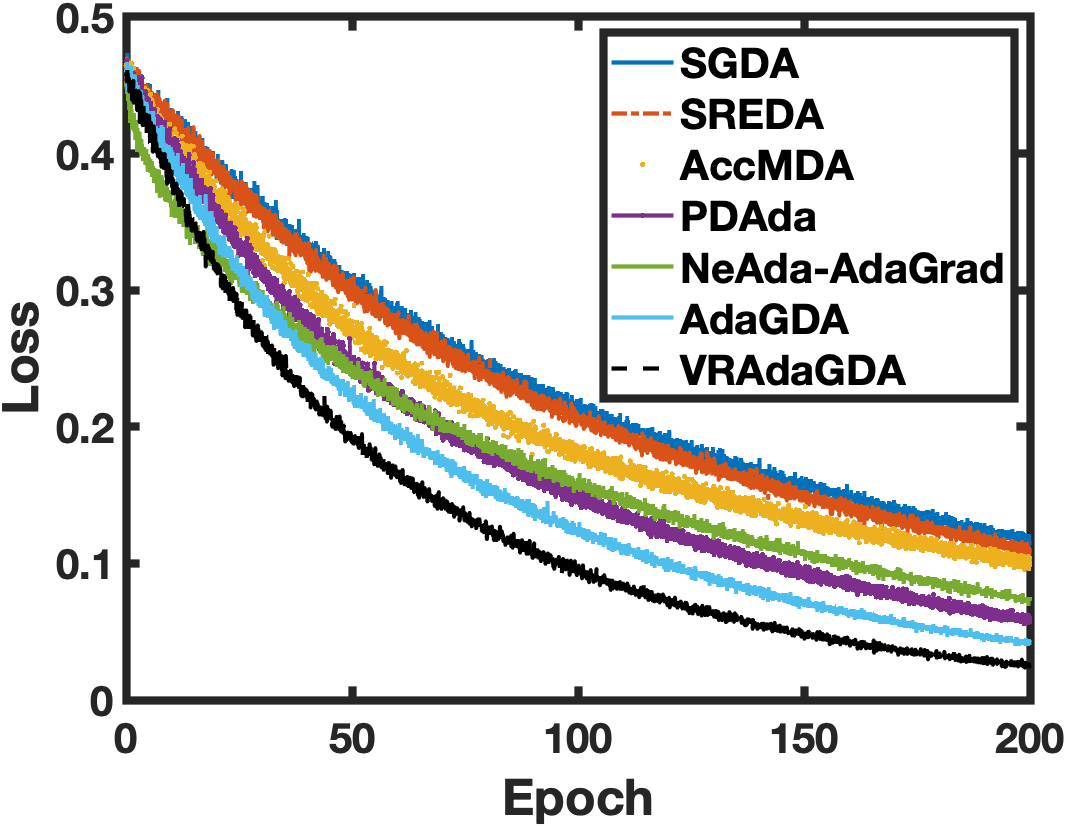}
}
\subfigure[MountainCarContinuous-v0]{
\hspace{0pt}
\includegraphics[width=.31\textwidth]{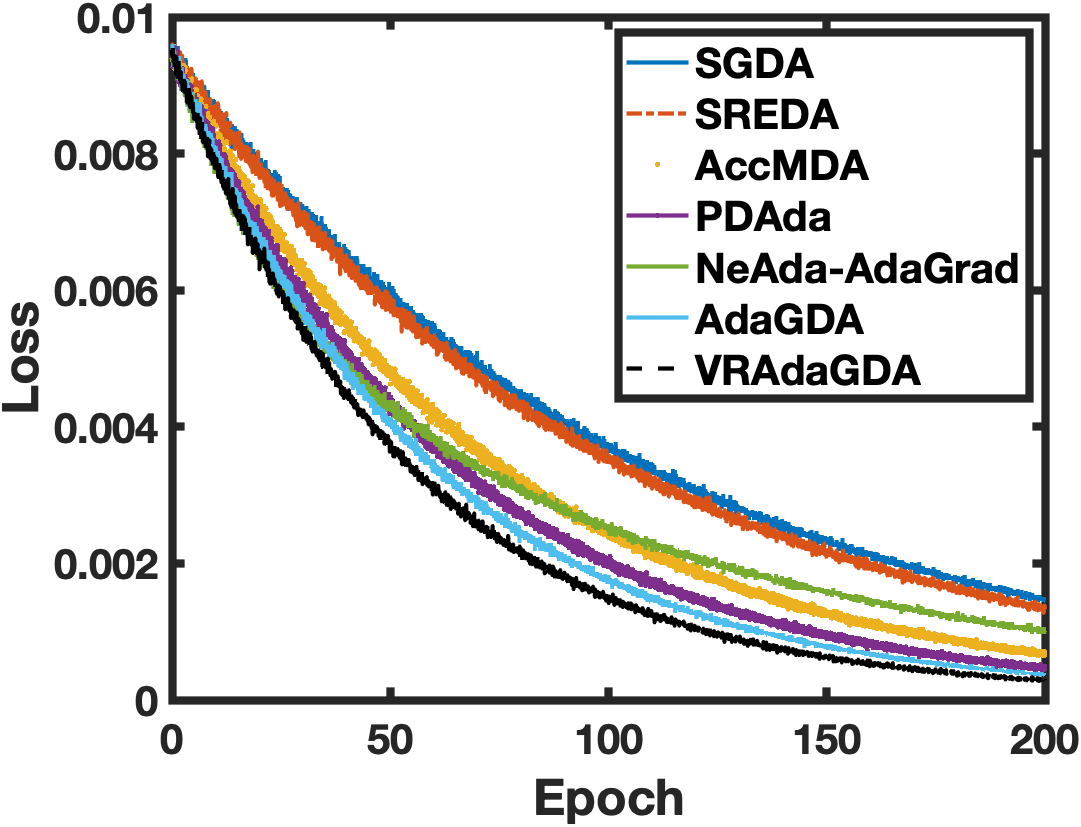}
}
\caption{Results of different methods on the policy evaluation task.}
\label{fig:PE}
\end{figure*}

\begin{figure*}[ht]
\centering
\subfigure[Fashion-MNIST]{
\hspace{0pt}
\includegraphics[width=.31\textwidth]{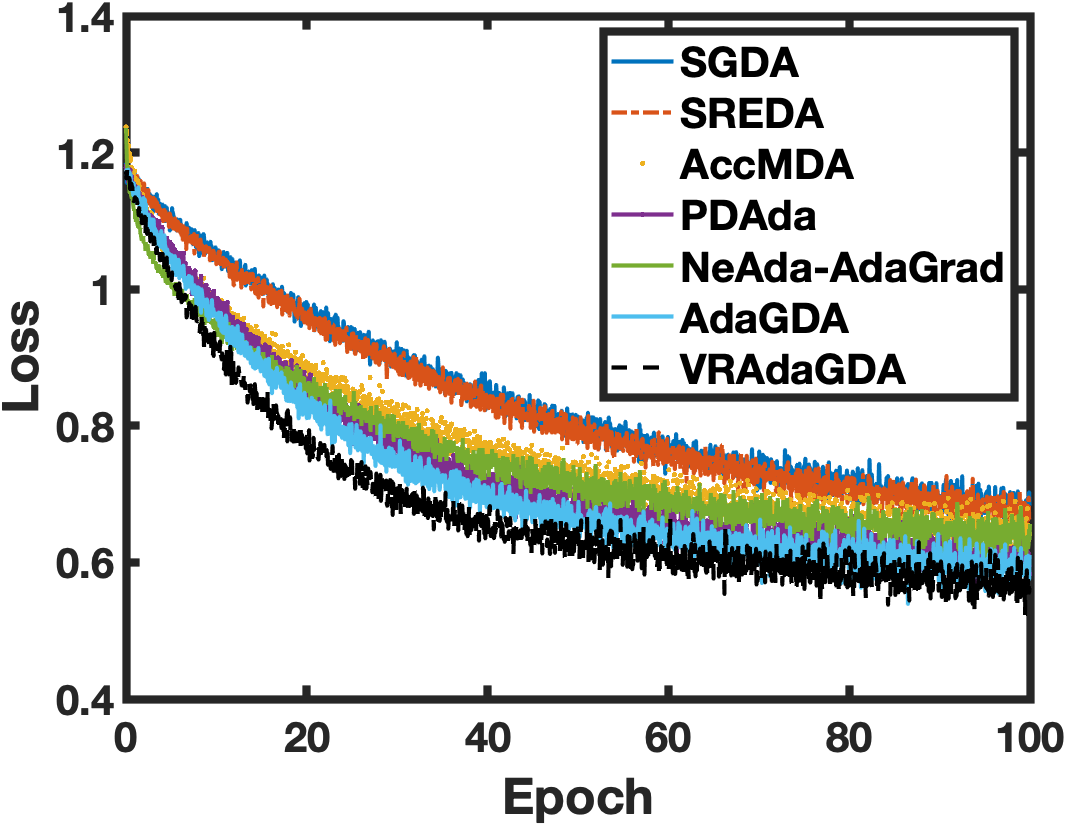}
}
\subfigure[MNIST]{
\hspace{0pt}
\includegraphics[width=.31\textwidth]{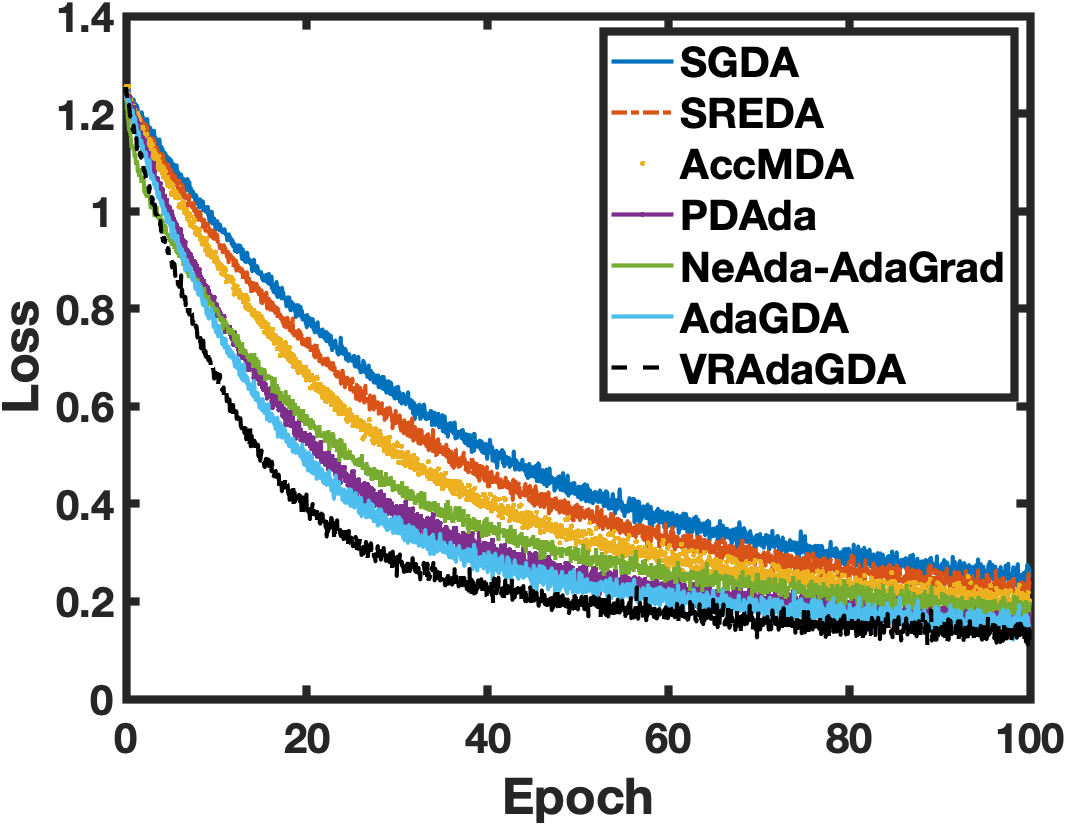}
}
\subfigure[CIFAR-10]{
\hspace{0pt}
\includegraphics[width=.31\textwidth]{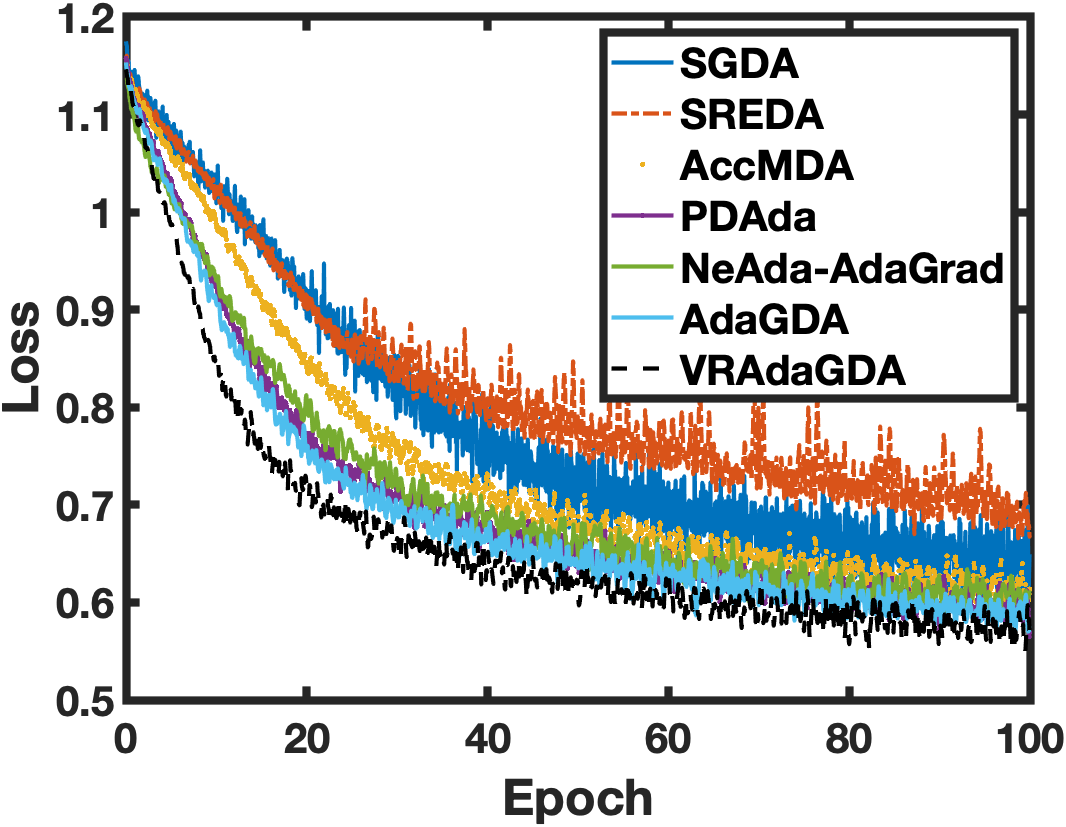}
}
\caption{Results of different methods on the fair classifier task.}
\label{fig:fairclass}
\end{figure*}

\begin{table}
  \centering
  \caption{ Model Architecture for the Policy Evaluation}
  \label{tab:2}
\begin{tabular}{ll}
\hline Layer Type & Shape \\
\hline
Fully Connected + tanh & 16 \\
Fully Connected       & 1 \\
\hline
\end{tabular}
\end{table}

\begin{table}
  \centering
  \caption{ Model Architecture for the Fair Classifier }
  \label{tab:3}
\begin{tabular}{ll}
\hline Layer Type & Shape \\
\hline Convolution + ReLU & $3 \times 3 \times 5 $ \\
Max Pooling & $2 \times 2$ \\
Convolution + ReLU & $3 \times 3 \times 10$ \\
Max Pooling & $2 \times 2$ \\
Fully Connected + ReLU & 100 \\
Fully Connected + ReLU & 3 \\
\hline
\end{tabular}
\end{table}

\subsection{ Policy Evaluation }
The first task is to apply a neural network to estimate the value function in Markov Decision Process (MDP). The value function $V_{\boldsymbol{\theta}}(\cdot) $
is parameterized as a 2-layer neural network, whose minimax loss function is defined in \eqref{eq:2} given in the Introduction.
In the experiment, we generate 10,000 state-reward pairs for three classic environments from GYM \citep{brockman2016openai}: CartPole-v1, Acrobat-v1, and MountainCarContinuous-v0.
Specifically, in CartPole-v1, a pole is connected with a cart by a joint. The goal of CartPole-v1 is to keep the pole upright by adding force to the cart. The system in  Acrobot-v1 has two joints and two links. To get the reward, we need to swing the end of the lower link and make it reach a given height.
In MountainCarContinuous-v0, the car is on a one-dimensional track between two "mountains". The car needs to drive up to the mountain on the right
but the car's engine is not strong enough to complete this task without momentum.

In the MDP, we let the discount factor $\tau=0.95$. In our algorithms, we set $\gamma=\lambda=0.005$, and the adaptive matrices $A_t$ and $B_t$
are generated from \eqref{eq:4} and \eqref{eq:5} respectively, where $\varrho=0.1$ and $\rho  =0.001$. In other algorithms,
we set the step-size for updating parameter $\theta$ be 0.005 and the step-size for $\omega$ be 0.005. At the same time, in the SREDA algorithm,
we set $S_1=10,000$ and $S_2=q=500$. The batch-sizes for all other methods are 500. In AccMDA and VR-AdaGDA, $\alpha_{t+1}= \eta^2_t$, $\beta_{t+1}= \eta^2_t$.
In  AdaGDA, $\alpha_{t+1}= \eta_t$, $\beta_{t+1}=\eta_t$. In PDAda, $\beta_{x} = \beta_{t} = \eta_x = \eta_y = 0.9$.
In NeAda-AdaGrad~\citep{yang2022nest}, we utilized the AdaGrad~\citep{duchi2011adaptive} optimizer in both dual and prime variables. The step-size is chosen from the set $0.015$. To avoid the explosion of adaptive learning rates, we clip it between $(0,3)$.
The architecture of neural network for policy evaluation is given in Table~\ref{tab:2}.

Figure \ref{fig:PE} shows the loss \emph{vs.} epoch of different stochastic methods. From these results, we can observe that our algorithms outperform the other algorithms,
and the VR-AdaGDA  consistently outperforms the AdaGDA.

\subsection{ Fair Classifier }
In the second task, we train a fair classifier by minimizing the maximum loss over different categories, where we use a Convolutional Neural Network (CNN) model as classifier.
In the experiment, we use the MNIST, Fashion-MNIST, and CIFAR-10 datasets as in \cite{nouiehed2019solving}. Following \cite{nouiehed2019solving},
we mainly focus on three categories in each dataset: digital numbers $\{0, 2, 3\}$ in the MNIST dataset, and T-shirt/top,
Coat and Shirt categories in the Fashion-MNIST dataset, and airplane, automobile and bird in the CIFAR10 dataset. Then we train this fair classifier by solving the following minimax problem:
\begin{align} \label{eq:23}
 \min_{w} \max_{ u\in \mathcal{U} } \ \big\{\sum_{i=1}^{3} u_i \mathcal{L}_i(w) - \varrho\|u-\frac{\textbf{1}}{3}\|^2  \big\},
\end{align}
where $\mathcal{U}=\big\{ u \ | \  u_i \geq 0, \ \sum_{i=1}^{3} u_{i}=1 \big\}$, $\mathcal{L}_{1}$, $\mathcal{L}_{2}$, and $\mathcal{L}_{3}$ are the cross-entropy loss functions corresponding to the samples in three different categories. Here $\varrho \geq 0$ is tuning parameter, and $u$ is a weight vector for different loss functions, and $w$ denotes the parameters of CNN.

In the experiment, we use xavier normal initialization to CNN  layer.
In our algorithms, we set $\gamma=0.001$ and $\lambda=0.0001$, and the adaptive matrices $A_t$ and $B_t$ are generated from \eqref{eq:4} and \eqref{eq:5} respectively, where $\varrho=0.1$ and $\rho =0.001$. In the other algorithms, we set the step-size for updating parameter $w$ be 0.001 and step-size for $u$ be  0.0001. At the same time, we set $\eta_t=0.9$ in our algorithms.
We run all algorithms for 100 epochs, and then record the loss value. For SREDA, we set $S_1=18,000$ and $S_2=q=900$. The batch-sizes for all other methods are 900. For AccMDA and VR-AdaGDA, $\alpha_{t+1}= \eta^2_t$, $\beta_{t+1}= \eta^2_t$. For AdaGDA, $\alpha_{t+1}= \eta_t$, $\beta_{t+1}=\eta_t$. For PDAda,  $\beta_{x} = \beta_{t} = \eta_x = \eta_y = 0.9$. In NeAda-AdaGrad, we utilized the AdaGrad optimizer in both dual and prime variables. The step-size is set as $0.015$.  Note that for fair comparison, we do not use the small stepsizes relying on small $\epsilon$ following the original SREDA algorithm, but use the relatively large stepsizes in the experiments. The architecture of CNN for policy evaluation is given in Table~\ref{tab:3}.

Figure \ref{fig:fairclass} plots the loss \emph{vs.} epoch of different stochastic methods. From these results, we can see that our algorithms consistently outperform other related methods.

\section{Conclusions}
In the paper, we proposed a class of faster adaptive gradient descent ascent methods for solving the minimax Problem \eqref{eq:1} using unified adaptive matrices
for both variables $x$ and $y$. In particular, our methods can easily incorporate both the momentum and variance-reduced techniques. Moreover, we provided an effective convergence analysis framework for our proposed methods, and proved that our methods obtain the best known gradient complexity for finding the first-order stationary points. The empirical studies on policy evaluation and fair classifier learning tasks were conducted to validate the efficiency of our new algorithms.

\subsubsection*{Acknowledgements}
We thank the anonymous reviewers for their valuable comments.
We also thank for the help of Prof. Heng Huang.
This work was partially supported by NSFC under Grant No.
61806093. Feihu Huang is the corresponding author (huangfeihu2018@gmail.com).

\small

\bibliography{AdaGDA}

\appendix
\onecolumn

\section{Appendix}
In this section, we provide the detailed convergence analysis of our algorithms.
We first give some useful lemmas.

Given a $\rho$-strongly convex function $\psi(x):
\mathcal{X}\rightarrow \mathbb{R}$, we define a Bregman distance \citep{censor1981iterative,censor1992proximal,ghadimi2016mini} associated with $\psi(x)$ as follows:
\begin{align}
 D(z,x) = \psi(z) - \big[\psi(x) + \langle\nabla \psi(x), z-x\rangle\big], \quad \forall{x,z} \in \mathcal{X},
\end{align}
where $\mathcal{X}\subseteq \mathbb{R}^d$ is a closed convex set.
Assume $h(x): \mathcal{X} \rightarrow \mathbb{R}$  is a convex and possibly nonsmooth function,
we define a generalized projection problem:
\begin{align} \label{eq:A1}
 x^+ = \arg\min_{z\in \mathcal{X}} \big\{\langle z,v\rangle + h(z) + \frac{1}{\gamma}D(z,x)\big\}, \quad x\in \mathcal{X},
\end{align}
where $v\in \mathbb{R}^d$ and $\gamma>0$.
Following \cite{ghadimi2016mini}, we define a generalized gradient as follows:
\begin{align}
 \mathcal{G}_{\mathcal{X}}(x,v,\gamma) = \frac{1}{\gamma}(x-x^+).
\end{align}

\begin{lemma} \label{lem:A1}
(Lemma 1 in \cite{ghadimi2016mini})
Let $x^+$ be given in \eqref{eq:A1}. Then we
have, for any $x\in \mathcal{X}$, $v\in \mathbb{R}^d$ and $\gamma >0$,
\begin{align}
 \langle v,  \mathcal{G}_{\mathcal{X}}(x,v,\gamma)\rangle \geq \rho \|\mathcal{G}_{\mathcal{X}}(x,v,\gamma)\|^2 + \frac{1}{\gamma}\big[h(x^+)-h(x)\big],
\end{align}
where $\rho>0$ depends on $\rho$-strongly convex function $\psi(x)$.
\end{lemma}
Based on Lemma \ref{lem:A1}, let $h(x)=0$, we have
\begin{align}
 \langle v,  \mathcal{G}_{\mathcal{X}}(x,v,\gamma)\rangle \geq \rho \|\mathcal{G}_{\mathcal{X}}(x,v,\gamma)\|^2.
\end{align}

\begin{lemma} \label{lem:A2}
\citep{nesterov2018lectures}
Assume function $f(x)$ is convex and $\mathcal{X}$ is a convex set.
 $x^* \in \mathcal{X}$ is the solution of the
constrained problem $\min_{x\in \mathcal{X}}f(x)$, if
\begin{align}
 \langle \nabla f(x^*), x-x^*\rangle \geq 0, \ \forall x\in \mathcal{X}.
\end{align}
where $\nabla f(x^*)$ denote the (sub-)gradient of function $f(x)$ at $x^*$.
\end{lemma}

\begin{lemma} \label{lem:A3}
\citep{lin2019gradient} Under the above Assumptions \ref{ass:2} and \ref{ass:5}, the function $F(x)=\min_{y\in \mathcal{Y}} f(x,y) =f(x,y^*(x))$ and the mapping $y^*(x)=\arg\max_{y\in \mathcal{Y}}f(x,y)$ have $L$-Lipschitz continuous gradient and $\kappa$-Lipschitz continuous respectively, such as for all $x_1,x_2 \in \mathcal{X}$
\begin{align}
 \|\nabla F(x_1) - \nabla F(x_2)\| \leq L\|x_1-x_2\|, \quad \|y^*(x_1) - y^*(x_2)\| \leq \kappa \|x_1-x_2\|,
\end{align}
where $L=L_f (1+ \kappa)$ and $\kappa = L_f/\mu$.
\end{lemma}

\begin{lemma} \label{lem:A4}
For independent random variables $\{\xi_i\}_{i=1}^n$ with zero mean, we have
$\mathbb{E}\|\frac{1}{n}\sum_{i=1}^n\xi_i\|^2 = \frac{1}{n}\mathbb{E}\|\xi_i\|^2$ for any $i\in [n]$.
\end{lemma}

\begin{lemma} \label{lem:D1}
Suppose the sequence $\{x_t,y_t\}_{t=1}^T$ be generated from Algorithms \ref{alg:1} or \ref{alg:2}. Let $0<\eta_t\leq 1$ and $0< \gamma \leq \frac{\rho}{2L\eta_t}$,
we have
\begin{align}
F(x_{t+1}) - F(x_t) \leq \frac{2\gamma L_f^2\eta_t}{\rho}\|y^*(x_t)-y_t\|^2 + \frac{2\gamma\eta_t}{\rho}\|\nabla_{x} f(x_t,y_t) -v_t \|^2 -\frac{\rho\eta_t}{2\gamma}\|\tilde{x}_{t+1}-x_t\|^2,
\end{align}
where $L = L_f (1+ \kappa)$.
\end{lemma}

\begin{proof}
According to the above Lemma \ref{lem:A3}, the function $F(x)$ has $L$-Lipschitz continuous gradient.
Then we have
\begin{align} \label{eq:D1}
  F(x_{t+1}) &\leq F(x_t) + \langle\nabla F(x_t), x_{t+1}-x_t\rangle + \frac{L}{2}\|x_{t+1}-x_t\|^2  \nonumber \\
  & = F(x_t) + \eta_t\langle \nabla F(x_t),\tilde{x}_{t+1}-x_t\rangle + \frac{L\eta_t^2}{2}\|\tilde{x}_{t+1}-x_t\|^2 \nonumber \\
  & = F(x_t) + \eta_t \underbrace{\langle v_t ,\tilde{x}_{t+1}-x_t\rangle}_{=T_1}+ \eta_t\underbrace{\langle \nabla F(x_t)-v_t ,\tilde{x}_{t+1}-x_t\rangle}_{=T_2} + \frac{L\eta_t^2}{2}\|\tilde{x}_{t+1}-x_t\|^2,
\end{align}
where the first equality holds by $x_{t+1}=x_t + \eta_t(\tilde{x}_{t+1}-x_t)$.

According to Assumption \ref{ass:4}, i.e., $A_t\succ \rho I_{d_1}$ for any $t\geq 1$,
the function $\phi_t(x)=x^TA_tx$ is $\rho$-strongly convex.
By using the above Lemma \ref{lem:A1} to the line 5 of Algorithm \ref{alg:1} or \ref{alg:2} , we have
\begin{align}
  \langle v_t , \frac{1}{\gamma}(x_t - \tilde{x}_{t+1})\rangle \geq \rho\|\frac{1}{\gamma}(x_t - \tilde{x}_{t+1})\|^2 \Rightarrow \langle v_t , \tilde{x}_{t+1}-x_t\rangle \leq -\frac{\rho}{\gamma}\|\tilde{x}_{t+1}-x_t\|^2.
\end{align}
Then we obtain
\begin{align} \label{eq:D2}
  T_1=\langle v_t , \tilde{x}_{t+1}-x_t\rangle \leq -\frac{\rho}{\gamma}\|\tilde{x}_{t+1}-x_t\|^2.
\end{align}

Next, we decompose the term $T_2=\langle \nabla F(x_t)-v_t ,\tilde{x}_{t+1}-x_t\rangle$ as follows:
\begin{align}
T_2 & =\langle \nabla F(x_t)-v_t ,\tilde{x}_{t+1}-x_t\rangle \nonumber \\
& = \underbrace{\langle \nabla F(x_t) - \nabla_{x} f(x_t,y_t),\tilde{x}_{t+1}-x_t\rangle}_{=T_3}
+ \underbrace{\langle  \nabla_{x} f(x_t,y_t) -v_t ,\tilde{x}_{t+1}-x_t\rangle}_{=T_4}.
\end{align}
For the term $T_3$, by the Cauchy-Schwarz inequality and Young's inequality, we have
\begin{align}
T_3 &= \langle \nabla F(x_t) - \nabla_{x} f(x_t,y_t),\tilde{x}_{t+1}-x_t\rangle \nonumber \\
& \leq \|\nabla F(x_t) - \nabla_{x} f(x_t,y_t)\|\cdot\|\tilde{x}_{t+1}-x_t\| \nonumber \\
& \leq \frac{2\gamma}{\rho}\|\nabla F(x_t) - \nabla_{x} f(x_t,y_t)\|^2 + \frac{\rho}{8\gamma}\|\tilde{x}_{t+1}-x_t\|^2 \nonumber \\
& = \frac{2\gamma}{\rho}\|\nabla_{x} f(x_t,y^*(x_t)) - \nabla_{x} f(x_t,y_t)\|^2 + \frac{\rho}{8\gamma}\|\tilde{x}_{t+1}-x_t\|^2 \nonumber \\
& \leq \frac{2\gamma L_f^2}{\rho}\|y^*(x_t)-y_t\|^2 + \frac{\rho}{8\gamma}\|\tilde{x}_{t+1}-x_t\|^2,
\end{align}
where the second inequality is due to the inequality $\langle a,b\rangle \leq \frac{\nu}{2}\|a\|^2 + \frac{1}{2\nu}\|b\|^2$
with $\nu=\frac{4\gamma}{\rho}$, and the last inequality holds by Assumption \ref{ass:5}.
For the term $T_2$, similarly, we have
\begin{align}
T_4 & = \langle  \nabla_{x} f(x_t,y_t) -v_t ,\tilde{x}_{t+1}-x_t \rangle \nonumber \\
& \leq \|\nabla_{x} f(x_t,y_t) -v_t \| \cdot \|\tilde{x}_{t+1}-x_t\| \nonumber \\
& \leq \frac{2\gamma}{\rho}\|\nabla_{x} f(x_t,y_t) -v_t \|^2 + \frac{\rho}{8\gamma}\|\tilde{x}_{t+1}-x_t\|^2.
\end{align}
Thus, we have
\begin{align} \label{eq:D4}
T_2 = \frac{2\gamma L_f^2}{\rho}\|y^*(x_t)-y_t\|^2 + \frac{2\gamma}{\rho}\|\nabla_{x} f(x_t,y_t) -v_t \|^2 + \frac{\rho}{4\gamma}\|\tilde{x}_{t+1}-x_t\|^2.
\end{align}
Finally, combining the inequalities \eqref{eq:D1}, \eqref{eq:D2} with \eqref{eq:D4}, we have
\begin{align}
 F(x_{t+1}) & \leq F(x_t) -\frac{\rho\eta_t}{\gamma}\|\tilde{x}_{t+1}-x_t\|^2 + \frac{2\gamma L_f^2\eta_t}{\rho}\|y^*(x_t)-y_t\|^2 + \frac{2\gamma\eta_t}{\rho}\|\nabla_{x} f(x_t,y_t) -v_t \|^2 \nonumber \\
 & \quad + \frac{\rho\eta_t}{4\gamma}\|\tilde{x}_{t+1}-x_t\|^2 + \frac{L\eta_t^2}{2}\|\tilde{x}_{t+1}-x_t\|^2\nonumber \\
 & \leq F(x_t) + \frac{2\gamma L_f^2\eta_t}{\rho}\|y^*(x_t)-y_t\|^2 + \frac{2\gamma\eta_t}{\rho}\|\nabla_{x} f(x_t,y_t) -v_t \|^2 -\frac{\rho\eta_t}{2\gamma}\|\tilde{x}_{t+1}-x_t\|^2,
\end{align}
where the last inequality is due to $0< \gamma \leq \frac{\rho}{2L\eta_t}$.

\end{proof}

\begin{lemma} \label{lem:E1}
Suppose the sequence $\{x_t,y_t\}_{t=1}^T$ be generated from Algorithm \ref{alg:1} or \ref{alg:2} .
Under the above Assumptions, given $B_t=b_tI_{d_2} \ (b_t\geq b>0)$ for all $t\geq 1$, $0< \eta_t\leq 1$
and $0<\lambda \leq \frac{b}{6L_f} \leq \frac{b_t}{6L_f}$, we have
\begin{align}
\|y_{t+1} - y^*(x_{t+1})\|^2 & \leq (1-\frac{\eta_t\mu\lambda}{4b_t})\|y_t -y^*(x_t)\|^2 -\frac{3\eta_t}{4} \|\tilde{y}_{t+1}-y_t\|^2 \nonumber \\
     & \quad + \frac{25\eta_t\lambda}{6\mu b_t} \|\nabla_y f(x_t,y_t)-w_t\|^2 + \frac{25\kappa^2\eta_tb_t}{6\mu\lambda}\|\tilde{x}_{t+1} - x_t\|^2,
\end{align}
where $\kappa = L_f/\mu$.
\end{lemma}
\begin{proof}
 According to Assumption \ref{ass:2}, i.e., the function $f(x,y)$ is $\mu$-strongly concave w.r.t $y$,
 we have
 \begin{align} \label{eq:E1}
  f(x_t,y) & \leq f(x_t,y_t) + \langle\nabla_y f(x_t,y_t), y-y_t\rangle - \frac{\mu}{2}\|y-y_t\|^2 \nonumber \\
  & = f(x_t,y_t) + \langle w_t, y-\tilde{y}_{t+1}\rangle + \langle\nabla_y f(x_t,y_t)-w_t, y-\tilde{y}_{t+1}\rangle \nonumber \\
  & \quad +\langle\nabla_y f(x_t,y_t), \tilde{y}_{t+1}-y_t\rangle- \frac{\mu}{2}\|y-y_t\|^2.
 \end{align}
 According to Assumption \ref{ass:5}, i.e., the function $f(x,y)$ is $L_f$-smooth, we have
 \begin{align} \label{eq:E2}
  -\frac{L_f}{2}\|\tilde{y}_{t+1}-y_t\|^2 \leq f(x_t,\tilde{y}_{t+1}) - f(x_t,y_{t}) -\langle\nabla_y f(x_t,y_t), \tilde{y}_{t+1}-y_t\rangle.
 \end{align}
 Summing up the about inequalities \eqref{eq:E1} with \eqref{eq:E2}, we have
 \begin{align} \label{eq:E3}
  f(x_t,y) & \leq f(x_t,\tilde{y}_{t+1}) + \langle w_t, y-\tilde{y}_{t+1}\rangle + \langle\nabla_y f(x_t,y_t)-w_t, y-\tilde{y}_{t+1}\rangle \nonumber \\
  & \quad - \frac{\mu}{2}\|y-y_t\|^2 + \frac{L_f}{2}\|\tilde{y}_{t+1}-y_t\|^2.
 \end{align}

 By the optimality of the line 6 of Algorithm \ref{alg:1} or \ref{alg:2} and $B_t=b_tI_{d_2}$, we have
  \begin{align}
  \langle -w_t + \frac{b_t}{\lambda}(\tilde{y}_{t+1} - y_t), y-\tilde{y}_{t+1} \rangle \geq 0,  \quad \forall y\in \mathcal{Y}
  \end{align}
 where the above inequality holds by Lemma \ref{lem:A2}.
 Then we obtain
 \begin{align} \label{eq:E4}
  \langle w_t, y-\tilde{y}_{t+1}\rangle & \leq \frac{1}{\lambda}\langle b_t(\tilde{y}_{t+1}- y_t), y-\tilde{y}_{t+1}\rangle  \nonumber \\
  & = \frac{1}{\lambda}\langle b_t(\tilde{y}_{t+1}- y_t), y_t-\tilde{y}_{t+1}\rangle
  + \frac{1}{\lambda}\langle b_t(\tilde{y}_{t+1}- y_t), y-y_t\rangle  \nonumber \\
  & = -\frac{b_t}{\lambda}\|\tilde{y}_{t+1}- y_t\|^2 + \frac{b_t}{\lambda}\langle \tilde{y}_{t+1}- y_t, y-y_t\rangle.
  \end{align}

 By plugging the inequalities \eqref{eq:E4} into \eqref{eq:E3}, we have
 \begin{align}
  f(x_t,y) & \leq f(x_t,\tilde{y}_{t+1}) + \frac{b_t}{\lambda}\langle \tilde{y}_{t+1}- y_t, y-y_t\rangle + \langle\nabla_y f(x_t,y_t)-w_t, y-\tilde{y}_{t+1}\rangle \nonumber \\
  & \quad -\frac{b_t}{\lambda}\|\tilde{y}_{t+1}- y_t\|^2- \frac{\mu}{2}\|y-y_t\|^2 + \frac{L_f}{2}\|\tilde{y}_{t+1}-y_t\|^2.
 \end{align}
 Let $y=y^*(x_t)$ and we obtain
 \begin{align}
  f(x_t,y^*(x_t)) & \leq f(x_t,\tilde{y}_{t+1}) + \frac{b_t}{\lambda}\langle \tilde{y}_{t+1}- y_t, y^*(x_t)-y_t\rangle + \langle\nabla_y f(x_t,y_t)-w_t, y^*(x_t)-\tilde{y}_{t+1}\rangle \nonumber \\
  & \quad -\frac{b_t}{\lambda}\|\tilde{y}_{t+1}- y_t\|^2- \frac{\mu}{2}\|y^*(x_t)-y_t\|^2 + \frac{L_f}{2}\|\tilde{y}_{t+1}-y_t\|^2.
 \end{align}
 Due to the concavity of $f(\cdot,y)$ and $y^*(x_t) =\arg\max_{y\in \mathcal{Y}}f(x_t,y)$,
 we have $f(x_t,y^*(x_t)) \geq f(x_t,\tilde{y}_{t+1})$.
 Thus, we obtain
 \begin{align} \label{eq:E5}
  0 & \leq  \frac{b_t}{\lambda}\langle \tilde{y}_{t+1}- y_t, y^*(x_t)-y_t\rangle + \langle\nabla_y f(x_t,y_t)-w_t, y^*(x_t)-\tilde{y}_{t+1}\rangle \nonumber \\
  & \quad -\frac{b_t}{\lambda}\|\tilde{y}_{t+1}- y_t\|^2- \frac{\mu}{2}\|y^*(x_t)-y_t\|^2 + \frac{L_f}{2}\|\tilde{y}_{t+1}-y_t\|^2.
 \end{align}

 By $y_{t+1} = y_t + \eta_t(\tilde{y}_{t+1}-y_t) $, we have
 \begin{align}
  \|y_{t+1}-y^*(x_t)\|^2 & = \|y_t + \eta_t(\tilde{y}_{t+1}-y_t) -y^*(x_t)\|^2 \nonumber \\
  & = \|y_t -y^*(x_t)\|^2 + 2\eta_t\langle \tilde{y}_{t+1}-y_t, y_t -y^*(x_t)\rangle + \eta_t^2\|\tilde{y}_{t+1}-y_t\|^2.
 \end{align}
 Then we obtain
 \begin{align} \label{eq:E6}
  \langle \tilde{y}_{t+1}-y_t, y^*(x_t) - y_t\rangle \leq \frac{1}{2\eta_t}\|y_t -y^*(x_t)\|^2 + \frac{\eta_t}{2}\|\tilde{y}_{t+1}-y_t\|^2 - \frac{1}{2\eta_t}\|y_{t+1}-y^*(x_t)\|^2.
 \end{align}
 Considering the upper bound of the term $\langle\nabla_y f(x_t,y_t)-w_t, y^*(x_t)-\tilde{y}_{t+1}\rangle$, we have
 \begin{align} \label{eq:E7}
  &\langle\nabla_y f(x_t,y_t)-w_t, y^*(x_t)-\tilde{y}_{t+1}\rangle \nonumber \\
  & = \langle\nabla_y f(x_t,y_t)-w_t, y^*(x_t)-y_t\rangle + \langle\nabla_y f(x_t,y_t)-w_t, y_t-\tilde{y}_{t+1}\rangle \nonumber \\
  & \leq \frac{1}{\mu} \|\nabla_y f(x_t,y_t)-w_t\|^2 + \frac{\mu}{4}\|y^*(x_t)-y_t\|^2 + \frac{1}{\mu} \|\nabla_y f(x_t,y_t)-w_t\|^2 + \frac{\mu}{4}\|y_t-\tilde{y}_{t+1}\|^2 \nonumber \\
  & = \frac{2}{\mu} \|\nabla_y f(x_t,y_t)-w_t\|^2 + \frac{\mu}{4}\|y^*(x_t)-y_t\|^2 + \frac{\mu}{4}\|y_t-\tilde{y}_{t+1}\|^2.
 \end{align}

 By plugging the inequalities \eqref{eq:E6} and \eqref{eq:E7} into \eqref{eq:E5},
 we obtain
 \begin{align}
  \frac{b_t}{2\eta_t\lambda}\|y_{t+1}-y^*(x_t)\|^2 & \leq (\frac{b_t}{2\eta_t\lambda}-\frac{\mu}{4})\|y_t -y^*(x_t)\|^2 + \big( \frac{\eta_tb_t}{2\lambda} - \frac{b_t}{\lambda} + \frac{\mu}{4} + \frac{L_f}{2}\big) \|\tilde{y}_{t+1}-y_t\|^2  \\
  & \quad + \frac{2}{\mu} \|\nabla_y f(x_t,y_t)-w_t\|^2 \nonumber \\
  & \leq ( \frac{b_t}{2\eta_t\lambda}-\frac{\mu}{4})\|y_t -y^*(x_t)\|^2 + (\frac{3L_f}{4} -\frac{b_t}{2\lambda}) \|\tilde{y}_{t+1}-y_t\|^2 + \frac{2}{\mu} \|\nabla_y f(x_t,y_t)-w_t\|^2 \nonumber \\
  & = ( \frac{b_t}{2\eta_t\lambda}-\frac{\mu}{4})\|y_t -y^*(x_t)\|^2 - \big( \frac{3b_t}{8\lambda} + \frac{b_t}{8\lambda} -\frac{3L_f}{4}\big) \|\tilde{y}_{t+1}-y_t\|^2 \nonumber \\
  &\quad + \frac{2}{\mu} \|\nabla_y f(x_t,y_t)-w_t\|^2 \nonumber \\
  & \leq  \big(\frac{b_t}{2\eta_t\lambda}-\frac{\mu}{4}\big)\|y_t -y^*(x_t)\|^2 - \frac{3b_t}{8\lambda} \|\tilde{y}_{t+1}-y_t\|^2 + \frac{2}{\mu}\|\nabla_y f(x_t,y_t)-w_t\|^2, \nonumber
 \end{align}
 where the second inequality holds by $L_f \geq \mu$ and $0< \eta_t\leq 1$, and the last inequality is due to
 $0< \lambda \leq \frac{b}{6L_f} \leq \frac{b_t}{6L_f}$ for all $t\geq 1$.
 It implies that
 \begin{align} \label{eq:E8}
 \|y_{t+1}-y^*(x_t)\|^2 \leq (1-\frac{\eta_t\mu\lambda}{2b_t})\|y_t -y^*(x_t)\|^2 - \frac{3\eta_t}{4} \|\tilde{y}_{t+1}-y_t\|^2
 + \frac{4\eta_t\lambda}{\mu b_t}\|\nabla_y f(x_t,y_t)-w_t\|^2.
 \end{align}

 Next, we decompose the term $\|y_{t+1} - y^*(x_{t+1})\|^2$ as follows:
 \begin{align} \label{eq:E9}
   \|y_{t+1} - y^*(x_{t+1})\|^2
  & = \|y_{t+1} - y^*(x_t) + y^*(x_t) - y^*(x_{t+1})\|^2    \nonumber \\
  & =  \|y_{t+1} - y^*(x_t)\|^2 + 2\langle y_{t+1} - y^*(x_t), y^*(x_t) - y^*(x_{t+1})\rangle  + \|y^*(x_t) - y^*(x_{t+1})\|^2  \nonumber \\
  & \leq (1+\frac{\eta_t\mu\lambda}{4b_t})\|y_{t+1} - y^*(x_t)\|^2  + (1+\frac{4b_t}{\eta_t\mu\lambda})\|y^*(x_t) - y^*(x_{t+1})\|^2 \nonumber \\
  & \leq (1+\frac{\eta_t\mu\lambda}{4b_t})\|y_{t+1} - y^*(x_t)\|^2  + (1+\frac{4b_t}{\eta_t\mu\lambda})\kappa^2\|x_t - x_{t+1}\|^2,
 \end{align}
 where the first inequality holds by Cauchy-Schwarz inequality and Young's inequality, and  the second inequality is due to
 Lemma \ref{lem:A3}, and  the last equality holds by $x_{t+1}=x_t + \eta_t(\tilde{x}_{t+1}-x_t)$.

 By combining the above inequalities \eqref{eq:E8} and \eqref{eq:E9}, we have
 \begin{align}
  \|y_{t+1} - y^*(x_{t+1})\|^2 & \leq (1+\frac{\eta_t\mu\lambda}{4b_t})( 1-\frac{\eta_t\mu\lambda}{2b_t})\|y_t -y^*(x_t)\|^2
  - (1+\frac{\eta_t\mu\lambda}{4b_t})\frac{3\eta_t}{4} \|\tilde{y}_{t+1}-y_t\|^2     \nonumber \\
  & \quad + (1+\frac{\eta_t\mu\lambda}{4b_t})\frac{4\eta_t\lambda}{\mu b_t}\|\nabla_y f(x_t,y_t)-w_t\|^2
  + (1+\frac{4b_t}{\eta_t\mu\lambda})\kappa^2\|x_t - x_{t+1}\|^2.
 \end{align}
 Since $0 < \eta_t \leq 1$, $0< \lambda \leq \frac{b_t}{6L_f}$ and $L_f\geq \mu$, we have $\lambda \leq \frac{b_t}{6L_f} \leq \frac{b_t}{6\mu}$
 and $\eta_t\leq 1\leq \frac{b_t}{6\mu \lambda}$. Then we obtain
 \begin{align}
  (1+\frac{\eta_t\mu\lambda}{4b_t})( 1-\frac{\eta_t\mu\lambda}{2b_t})&= 1-\frac{\eta_t\mu\lambda}{2b_t} +\frac{\eta_t\mu\lambda}{4b_t}
  - \frac{\eta_t^2\mu^2\lambda^2}{8b_t^2} \leq 1-\frac{\eta_t\mu\lambda}{4b_t},  \\
  - (1+\frac{\eta_t\mu\lambda}{4b_t})\frac{3\eta_t}{4} &\leq -\frac{3\eta_t}{4}, \\
  (1+\frac{\eta_t\mu\lambda}{4b_t})\frac{4\eta_t\lambda}{\mu b_t} & \leq (1+\frac{1}{24})\frac{4\eta_t\lambda}{\mu}=\frac{25\eta_t\lambda}{6\mu b_t}, \\
  (1+\frac{4b_t}{\eta_t\mu\lambda})\kappa^2 &
  \leq \frac{\kappa^2b_t}{6\eta_t\mu\lambda} +\frac{4\kappa^2b_t}{\eta_t\mu\lambda} = \frac{25\kappa^2b_t}{6\eta_t\mu\lambda},
 \end{align}
 where the second last inequality is due to $\frac{\eta_t\mu\lambda}{b_t}\leq \frac{1}{6}$ and the last inequality holds by
 $\frac{b_t}{6\mu\lambda\eta_t} \geq 1$.
 Thus, we have
 \begin{align}
      \|y_{t+1} - y^*(x_{t+1})\|^2 &\leq (1-\frac{\eta_t\mu\lambda}{4b_t})\|y_t -y^*(x_t)\|^2 -\frac{3\eta_t}{4} \|\tilde{y}_{t+1}-y_t\|^2 \nonumber \\
      & \quad + \frac{25\eta_t\lambda}{6\mu b_t} \|\nabla_y f(x_t,y_t)-w_t\|^2 + \frac{25\kappa^2b_t}{6\mu\lambda \eta_t}\|x_{t+1} - x_t\|^2 \nonumber \\
      & = (1-\frac{\eta_t\mu\lambda}{4b_t})\|y_t -y^*(x_t)\|^2 -\frac{3\eta_t}{4} \|\tilde{y}_{t+1}-y_t\|^2 \nonumber \\
      & \quad + \frac{25\eta_t\lambda}{6\mu b_t} \|\nabla_y f(x_t,y_t)-w_t\|^2 + \frac{25\kappa^2\eta_tb_t}{6\mu\lambda}\|\tilde{x}_{t+1} - x_t\|^2,
 \end{align}
 where the equality holds by $x_{t+1}=x_t + \eta_t(\tilde{x}_{t+1}-x_t)$.

\end{proof}

\subsection{  Convergence Analysis of the AdaGDA Algorithm }
\label{Appendix:A1}
In this subsection, we study the convergence properties of our AdaGDA algorithm for solving the minimax problem \eqref{eq:1}.
We first give a useful Lemma for the gradient estimators.

\begin{lemma} \label{lem:C1}
 Assume that the stochastic partial derivatives $v_{t+1}$ and $w_{t+1}$ be generated from Algorithm \ref{alg:1}, we have
 \begin{align}
 \mathbb{E}\|\nabla_x f(x_{t+1},y_{t+1}) - v_{t+1}\|^2
 & \leq (1-\alpha_{t+1}) \mathbb{E} \|\nabla_x f(x_t,y_t) - v_t\|^2 + \frac{\alpha_{t+1}^2\sigma^2}{q} \nonumber \\
 & \quad + \frac{2L_f^2\eta_t^2}{\alpha_{t+1}}\big(\mathbb{E}\|\tilde{x}_{t+1} - x_t\|^2 + \mathbb{E}\|\tilde{y}_{t+1} - y_t\|^2 \big), \nonumber
 \end{align}
 \begin{align}
 \mathbb{E}\|\nabla_y f(x_{t+1},y_{t+1}) - w_{t+1}\|^2
 & \leq (1-\beta_{t+1}) \mathbb{E}\|\nabla_{y} f(x_t,y_t) -w_t\|^2 + \frac{\beta_{t+1}^2\sigma^2}{q} \nonumber \\
 & \quad + \frac{2L^2_f\eta^2_t}{\beta_{t+1}}\big(\mathbb{E}\|\tilde{x}_{t+1} - x_t\|^2 + \mathbb{E}\|\tilde{y}_{t+1}-y_t\|^2\big). \nonumber
 \end{align}
\end{lemma}

\begin{proof}
We first consider the term $\mathbb{E}\|\nabla_x f(x_{t+1},y_{t+1})-v_{t+1}\|^2$.
Since $v_{t+1} =  \alpha_{t+1}\nabla_x f(x_{t+1},y_{t+1};\mathcal{B}_{t+1}) + (1-\alpha_{t+1})v_t$, we have
\begin{align}
 &\mathbb{E}\|\nabla_x f(x_{t+1},y_{t+1})-v_{t+1}\|^2 \\
 & =  \mathbb{E}\|\nabla_x f(x_{t+1},y_{t+1}) -  \alpha_{t+1}\nabla_x f(x_{t+1},y_{t+1};\mathcal{B}_{t+1}) - (1-\alpha_{t+1})v_t \|^2 \nonumber \\
 & = \mathbb{E}\|\alpha_{t+1}(\nabla_x f(x_{t+1},y_{t+1}) - \nabla_x f(x_{t+1},y_{t+1};\mathcal{B}_{t+1})) + (1-\alpha_{t+1})(\nabla_x f(x_t,y_t) - v_t)\nonumber \\
 & \quad + (1-\alpha_{t+1})\big( \nabla_x f(x_{t+1},y_{t+1}) - \nabla_x f(x_t,y_t) \big)\|^2 \nonumber \\
 & = \mathbb{E}\| (1-\alpha_{t+1})(\nabla_x f(x_t,y_t) - v_t) + (1-\alpha_{t+1})\big( \nabla_x f(x_{t+1},y_{t+1}) - \nabla_x f(x_t,y_t) \big)\|^2 \nonumber \\
 & \quad + \alpha^2_{t+1}\mathbb{E}\|\nabla_x f(x_{t+1},y_{t+1}) - \nabla_x f(x_{t+1},y_{t+1};\mathcal{B}_{t+1})\|^2 \nonumber \\
 & \leq  (1-\alpha_{t+1})^2(1+\frac{1}{\alpha_{t+1}})\mathbb{E} \|\nabla_x f(x_{t+1},y_{t+1}) - \nabla_x f(x_t,y_t)\|^2 \nonumber \\
 & \quad +(1-\alpha_{t+1})^2(1+\alpha_{t+1})\mathbb{E} \|\nabla_x f(x_t,y_t) - v_t\|^2+ \alpha^2_{t+1}\mathbb{E}\|\nabla_x f(x_{t+1},y_{t+1}) - \nabla_x f(x_{t+1},y_{t+1};\mathcal{B}_{t+1})\|^2 \nonumber \\
 & \leq (1-\alpha_{t+1})\mathbb{E} \|\nabla_x f(x_t,y_t) - v_t\|^2 + \frac{1}{\alpha_{t+1}}\mathbb{E}\|\nabla_x f(x_{t+1},y_{t+1}) - \nabla_x f(x_t,y_t)\|^2  + \frac{\alpha^2_{t+1}\sigma^2}{q} \nonumber \\
 & \leq (1-\alpha_{t+1}) \mathbb{E} \|\nabla_x f(x_t,y_t) - v_t\|^2 + \frac{2L_f^2\eta_t^2}{\alpha_{t+1}}\big(\mathbb{E}\|\tilde{x}_{t+1} - x_t\|^2 + \mathbb{E}\|\tilde{y}_{t+1} - y_t\|^2 \big) + \frac{\alpha_{t+1}^2\sigma^2}{q}, \nonumber
\end{align}
where the third equality is due to $\mathbb{E}_{\mathcal{B}_{t+1}}[\nabla f(x_{t+1},y_{t+1};\mathcal{B}_{t+1})]=\nabla f(x_{t+1},y_{t+1})$; the second last inequality holds by $0\leq \alpha_{t+1} \leq 1$ such that  $(1-\alpha_{t+1})^2(1+\alpha_{t+1})=1-\alpha_{t+1}-\alpha_{t+1}^2+
  \alpha_{t+1}^3\leq 1-\alpha_{t+1}$ and $(1-\alpha_{t+1})^2(1+\frac{1}{\alpha_{t+1}}) \leq (1-\alpha_{t+1})(1+\frac{1}{\alpha_{t+1}})  =-\alpha_{t+1}+\frac{1}{\alpha_{t+1}}\leq \frac{1}{\alpha_{t+1}}$, and the last inequality holds by Assumption \ref{ass:5} and $x_{t+1}=x_t-\eta_t(\tilde{x}_{t+1}-x_t)$, $y_{t+1}=y_t-\eta_t(\tilde{y}_{t+1}-y_t)$.

Similarly, we have
\begin{align}
 \mathbb{E} \|\nabla_{y} f(x_{t+1},y_{t+1}) - w_{t+1}\|^2 & \leq (1-\beta_{t+1}) \mathbb{E}\|\nabla_{y} f(x_t,y_t) -w_t\|^2 + \frac{\beta_{t+1}^2\sigma^2}{q} \nonumber \\
 & \quad + \frac{2L^2_f\eta^2_t}{\beta_{t+1}}\big(\mathbb{E}\|\tilde{x}_{t+1} - x_t\|^2 + \mathbb{E}\|\tilde{y}_{t+1}-y_t\|^2\big).
\end{align}

\end{proof}

\begin{theorem} \label{th:A1}
(Restatement of Theorem 1)
Assume that the sequence $\{x_t,y_t\}_{t=1}^T$ be generated from the Algorithm \ref{alg:1}. When $\mathcal{X} \subset \mathbb{R}^{d_1}$, and given $B_t=b_tI_{d_2} \ (\hat{b} \geq b_t \geq b>0)$ for all $t\geq 1$,  $\eta_t=\frac{k}{(m+t)^{1/2}}$ for all $t\geq 0$, $\alpha_{t+1}=c_1\eta_t$, $\beta_{t+1}=c_2\eta_t$, $m\geq \max\big(k^2, (c_1k)^2, (c_2k)^2\big)$, $k>0$, $\frac{9\mu^2}{4} \leq c_1 \leq \frac{m^{1/2}}{k}$, $\frac{75L^2_f}{2} \leq c_2 \leq \frac{m^{1/2}}{k}$, $0< \gamma \leq \min\big(\frac{15\sqrt{2}\lambda\mu^2\rho}{2\sqrt{400L^2_f\lambda^2+24\mu^2\lambda^2+16875\hat{b}^2\kappa^2L^2_f\mu^2}}, \frac{m^{1/2}\rho}{4Lk}\big)$ and $0<\lambda\leq \min\big(\frac{405bL^2_f\mu^{3/2}}{8\sqrt{50L^2_f+9\mu^2}},\frac{b}{6L_f}\big)$, we have
\begin{align}
 \frac{1}{T} \sum_{t=1}^T\mathbb{E}\|\mathcal{G}_{\mathcal{X}}(x_t,\nabla F(x_t),\gamma)\|  \leq \frac{2\sqrt{3G}m^{1/4}}{T^{1/2}} + \frac{2\sqrt{3G}}{T^{1/4}},
\end{align}
where $G = \frac{F(x_1) - F^*}{k\gamma\rho} + \frac{9b_1L^2_f\Delta^2_1}{k\lambda\mu\rho^2} + \frac{2\sigma^2}{q k\mu^2\rho^2} + \frac{2m\sigma^2}{q k\mu^2\rho^2}\ln(m+T)$ and $\Delta^2_1=\|y_1-y^*(x_1)\|^2$.
\end{theorem}

\begin{proof}
Since $\eta_t=\frac{k}{(m+t)^{1/2}}$ on $t$ is decreasing and $m\geq k^2$, we have $\eta_t \leq \eta_0 = \frac{k}{m^{1/2}} \leq 1$ and $\gamma \leq \frac{m^{1/2}\rho}{4Lk}\leq \frac{\rho}{2L\eta_0} \leq \frac{\rho}{2L\eta_t}$ for any $t\geq 0$.
 Due to $0 < \eta_t \leq 1$ and $m\geq (c_1k)^2$, we have $\alpha_{t+1} = c_1\eta_t \leq \frac{c_1k}{m^{1/2}}\leq 1$.
 Similarly, due to $m\geq (c_2k)^2 $, we have $\beta_{t+1}\leq 1$. At the same time, we have $c_1,c_2\leq \frac{m^{1/2}}{k}$.
According to Lemma \ref{lem:C1}, we have
 \begin{align} \label{eq:H1}
  & \mathbb{E} \|\nabla_x f(x_{t+1},y_{t+1}) - v_{t+1}\|^2 - \mathbb{E} \|\nabla_x f(x_t,y_t) - v_t\|^2  \\
  & \leq -\alpha_{t+1}\mathbb{E} \|\nabla_x f(x_t,y_t) -v_t\|^2 + 2L^2_f\eta^2_t/\alpha_{t+1}\mathbb{E}\big(\|\tilde{x}_{t+1}-x_t\|^2 + \|\tilde{y}_{t+1}-y_t\|^2\big) + \frac{\alpha_{t+1}^2\sigma^2}{q}  \nonumber \\
  & = -c_1 \eta_t\mathbb{E} \|\nabla_x f(x_t,y_t) -v_t\|^2 + 2L^2_f\eta_t/c_1\mathbb{E}\big(\|\tilde{x}_{t+1}-x_t\|^2 + \|\tilde{y}_{t+1}-y_t\|^2\big) + \frac{c_1^2\eta_t^2\sigma^2}{q} \nonumber \\
  & \leq -\frac{9\mu^2\eta_t}{4}\mathbb{E} \|\nabla_x f(x_t,y_t) -u_t\|^2 + \frac{8L^2_f\eta_t}{9\mu^2}\mathbb{E}\big(\|\tilde{x}_{t+1}-x_t\|^2 + \|\tilde{y}_{t+1}-y_t\|^2\big) + \frac{m\eta^2_t\sigma^2}{k^2q}, \nonumber
 \end{align}
 where the above equality holds by $\alpha_{t+1}=c_1\eta_t$, and the last inequality is due to $\frac{9\mu^2}{4} \leq c_1 \leq \frac{m^{1/2}}{k}$.
 Similarly, given $\frac{75L^2_f}{2} \leq c_2 \leq \frac{m^{1/2}}{k}$, we have
 \begin{align} \label{eq:H2}
  & \mathbb{E} \|\nabla_y f(x_{t+1},y_{t+1}) - w_{t+1}\|^2 -  \mathbb{E} \|\nabla_y f(x_t,y_t) - w_t\|^2  \\
  & \leq - \frac{75L^2_f\eta_t}{2}\mathbb{E} \|\nabla_y f(x_t,y_t) - w_t\|^2 +
  \frac{4\eta_t}{75}\mathbb{E}\big(\|\tilde{x}_{t+1}-x_t\|^2 + \|\tilde{y}_{t+1}-y_t\|^2\big) + \frac{m\eta^2_t\sigma^2}{k^2q}. \nonumber
 \end{align}

According to Lemma \ref{lem:D1}, we have
\begin{align} \label{eq:H3}
F(x_{t+1}) - F(x_t) \leq \frac{2\gamma L_f^2\eta_t}{\rho}\|y^*(x_t)-y_t\|^2 + \frac{2\gamma\eta_t}{\rho}\|\nabla_{x} f(x_t,y_t) -v_t \|^2 -\frac{\rho\eta_t}{2\gamma}\|\tilde{x}_{t+1}-x_t\|^2.
\end{align}
According to Lemma \ref{lem:E1}, we have
\begin{align} \label{eq:H4}
\|y_{t+1} - y^*(x_{t+1})\|^2 - \|y_t -y^*(x_t)\|^2  &
\leq -\frac{\eta_t\mu\lambda}{4b_t} \|y_t -y^*(x_t)\|^2 -\frac{3\eta_t}{4} \|\tilde{y}_{t+1}-y_t\|^2  \\
& \quad + \frac{25\eta_t\lambda}{6\mu b_t} \|\nabla_y f(x_t,y_t)-w_t\|^2 + \frac{25\kappa^2\eta_tb_t}{6\mu\lambda}\|\tilde{x}_{t+1} - x_t\|^2. \nonumber
\end{align}

Next, we define a \emph{Lyapunov} function, for any $t\geq 1$
 \begin{align}
 \Omega_t & = \mathbb{E}\big [F(x_t) + \frac{9b_tL^2_f\gamma}{\lambda\mu\rho}\|y_t-y^*(x_t)\|^2 + \frac{\gamma}{\rho\mu^2} \big(\|\nabla_x f(x_t,y_t)-v_t\|^2 + \|\nabla_y f(x_t,y_t)-w_t\|^2\big) \big]. \nonumber
 \end{align}
Then we have
 \begin{align}
 & \Omega_{t+1} - \Omega_t \nonumber \\
 & = \mathbb{E}\big[F(x_{t+1}) - F(x_t)\big] + \frac{9b_tL^2_f\gamma}{\lambda\mu\rho} \big( \mathbb{E}\|y_{t+1}-y^*(x_{t+1})\|^2 - \mathbb{E}\|y_t-y^*(x_t)\|^2 \big)
 + \frac{\gamma}{\rho\mu^2} \big( \mathbb{E}\|\nabla_x f(x_{t+1},y_{t+1})-v_{t+1}\|^2 \nonumber \\
 & \quad - \mathbb{E}\|\nabla_x f(x_t,y_t)-v_t\|^2 + \mathbb{E}\|\nabla_y f(x_{t+1},y_{t+1})-w_{t+1}\|^2
 - \mathbb{E}\|\nabla_y f(x_t,y_t)-w_t\|^2 \big) \nonumber \\
 & \leq \frac{2\gamma L_f^2\eta_t}{\rho}\mathbb{E}\|y^*(x_t)-y_t\|^2 + \frac{2\gamma\eta_t}{\rho}\mathbb{E}\|\nabla_{x} f(x_t,y_t) -v_t \|^2 -\frac{\rho\eta_t}{2\gamma}\mathbb{E}\|\tilde{x}_{t+1}-x_t\|^2 \nonumber \\
 & \quad + \frac{9b_tL^2_f\gamma}{\lambda\mu\rho} \bigg( - \frac{\eta_t\mu\lambda}{4b_t} \mathbb{E}\|y_t -y^*(x_t)\|^2 -\frac{3\eta_t}{4} \mathbb{E}\|\tilde{y}_{t+1}-y_t\|^2 + \frac{25\eta_t\lambda}{6\mu b_t} \mathbb{E}\|\nabla_y f(x_t,y_t)-w_t\|^2  \nonumber \\
 & \quad  + \frac{25\kappa^2\eta_tb_t}{6\mu\lambda}\mathbb{E}\|\tilde{x}_{t+1} - x_t\|^2 \bigg) + \frac{\gamma}{\rho\mu^2} \bigg( -\frac{9\mu^2\eta_t}{4}\mathbb{E} \|\nabla_x f(x_t,y_t) -v_t\|^2 + \frac{8L^2_f\eta_t}{9\mu^2}\mathbb{E}\big(\|\tilde{x}_{t+1}-x_t\|^2 + \|\tilde{y}_{t+1}-y_t\|^2\big)  \nonumber \\
 & \quad + \frac{m\eta^2_t\sigma^2}{k^2q} - \frac{75L^2_f\eta_t}{2}\mathbb{E} \|\nabla_y f(x_t,y_t) - w_t\|^2 +
  \frac{4\eta_t}{75}\mathbb{E}\big(\|\tilde{x}_{t+1}-x_t\|^2 + \|\tilde{y}_{t+1}-y_t\|^2\big) + \frac{m\eta^2_t\sigma^2}{k^2q} \bigg)\nonumber \\
 & = - \frac{\gamma\eta_t}{4\rho} \big( L^2_f \mathbb{E}\|y_t - y^*(x_t)\|^2 + \mathbb{E} \|\nabla_x f(x_t,y_t) -v_t\|^2\big) - \big( \frac{\rho}{2\gamma} - \frac{8L^2_f\gamma}{9\rho\mu^4} - \frac{4\gamma}{75\rho\mu^2} - \frac{225b^2_t\kappa^2L^2_f\gamma}{6\mu^2\lambda^2\rho}\big)\eta_t\mathbb{E}\|\tilde{x}_{t+1}-x_t\|^2  \nonumber \\
 & \quad - \big(\frac{27b_tL^2_f\gamma}{4\lambda\mu\rho} -
  \frac{8L^2_f\gamma}{9\rho\mu^4} - \frac{4\gamma}{75\rho\mu^2}\big)\eta_t\mathbb{E}\|\tilde{y}_{t+1}-y_t\|^2
  + \frac{2m\gamma\sigma^2}{k^2\rho\mu^2q}\eta^2_t \nonumber \\
 & \leq - \frac{\gamma\eta_t}{4\rho} \big( L^2_f\mathbb{E}\|y_t - y^*(x_t)\|^2 + \mathbb{E} \|\nabla_x f(x_t,y_t) -v_t\|^2 \big) - \frac{\rho\eta_t}{4\gamma}\mathbb{E}\|\tilde{x}_{t+1}-x_t\|^2 + \frac{2m\gamma\sigma^2}{k^2\rho\mu^2q}\eta^2_t,
 \end{align}
 where the first inequality holds by the above inequalities \eqref{eq:H1}, \eqref{eq:H2}, \eqref{eq:H3} and \eqref{eq:H4};
 the last inequality is due to $0< \gamma \leq \frac{15\sqrt{2}\lambda\mu^2\rho}{2\sqrt{400L^2_f\lambda^2+24\mu^2\lambda^2+16875\hat{b}^2\kappa^2L^2_f\mu^2}} \leq \frac{15\sqrt{2}\lambda\mu^2\rho}{2\sqrt{400L^2_f\lambda^2+24\mu^2\lambda^2+16875b_t^2\kappa^2L^2_f\mu^2}}$ and $0<\lambda\leq \frac{405bL^2_f\mu^{3/2}}{8\sqrt{50L^2_f+9\mu^2}} \leq \frac{405b_tL^2_f\mu^{3/2}}{8\sqrt{50L^2_f+9\mu^2}}$ for all $t\geq 1$.
Then we have
\begin{align} \label{eq:H5}
  \frac{L^2_f\eta_t}{4}\mathbb{E}\|y_t - y^*(x_t)\|^2 + \frac{\eta_t}{4}\mathbb{E} \|\nabla_x f(x_t,y_t) -v_t\|^2 + \frac{\rho^2\eta_t}{4\gamma^2}\mathbb{E}\|\tilde{x}_{t+1}-x_t\|^2 \leq \frac{\rho(\Omega_t - \Omega_{t+1})}{\gamma} + \frac{2m\sigma^2}{k^2\mu^2q}\eta^2_t.
\end{align}
Taking average over $t=1,2,\cdots,T$ on both sides of \eqref{eq:H5}, we have
\begin{align}
 & \frac{1}{T} \sum_{t=1}^T \mathbb{E} \big[  \frac{L^2_f\eta_t}{4}\|y_t - y^*(x_t)\|^2 + \frac{\eta_t}{4}\|\nabla_x f(x_t,y_t) -v_t\|^2 + \frac{\rho^2\eta_t}{4\gamma^2}\|\tilde{x}_{t+1}-x_t\|^2 \big]  \nonumber \\
 & \leq  \sum_{t=1}^T \frac{\rho(\Omega_t - \Omega_{t+1})}{T\gamma} + \frac{1}{T}\sum_{t=1}^T\frac{2m\sigma^2}{k^2\mu^2q}\eta^2_t.
\end{align}
Given $x_1\in \mathcal{X}$, $y_1\in \mathcal{Y}$ and $\Delta^2_1 = \|y_1-y^*(x_1)\|^2$, we have
\begin{align} \label{eq:H6}
 \Omega_1 &= F(x_1) + \frac{9b_1L^2_f\gamma}{\lambda\mu\rho}\|y_1-y^*(x_1)\|^2 + \frac{\gamma}{\rho\mu^2} \big(\mathbb{E}\|\nabla_x f(x_1,y_1)-v_1\|^2 + \mathbb{E}\|\nabla_y f(x_1,y_1)-w_1\|^2 \big) \nonumber \\
 & \leq F(x_1) + \frac{9b_1L^2_f\gamma\Delta^2_1}{\lambda\mu\rho} + \frac{2\gamma\sigma^2}{q\rho\mu^2},
\end{align}
where the above inequality holds by Assumption \ref{ass:1}.

Since $\eta_t$ is decreasing on $t$, i.e., $\eta_T^{-1} \geq \eta_t^{-1}$ for any $0\leq t\leq T$, we have
 \begin{align}
 & \frac{1}{T} \sum_{t=1}^T \mathbb{E}\big[ \frac{L^2_f}{4}\|y_t - y^*(x_t)\|^2 + \frac{1}{4}\|\nabla_x f(x_t,y_t) -v_t\|^2 + \frac{\rho^2}{4\gamma^2}\|\tilde{x}_{t+1}-x_t\|^2 \big] \nonumber \\
 & \leq  \frac{\rho}{T\gamma\eta_T} \sum_{t=1}^T\big(\Omega_t - \Omega_{t+1}\big) + \frac{1}{T\eta_T}\sum_{t=1}^T\frac{2m\sigma^2}{k^2\mu^2q}\eta^2_t \nonumber \\
 & \leq \frac{\rho}{T\gamma\eta_T} \big( F(x_1) + \frac{9b_1L^2_f\gamma\Delta^2_1}{\lambda\mu\rho} + \frac{2\gamma\sigma^2}{q\rho\mu^2} - F^* \big) + \frac{1}{T\eta_T}\sum_{t=1}^T\frac{2m\sigma^2}{k^2\mu^2q}\eta^2_t  \nonumber \\
 & \leq \frac{\rho(F(x_1) - F^*)}{T\gamma\eta_T} + \frac{9b_1L^2_f\Delta^2_1}{\lambda\mu\eta_TT} + \frac{2\sigma^2}{q\mu^2\eta_TT} + \frac{2m\sigma^2}{\eta_TTk^2\mu^2q}\int^T_1\frac{k^2}{m+t} dt\nonumber \\
 & \leq  \frac{\rho(F(x_1) - F^*)}{T\gamma\eta_T} + \frac{9b_1L^2_f\Delta^2_1}{\lambda\mu\eta_TT} + \frac{2\sigma^2}{q\mu^2\eta_TT} + \frac{2m\sigma^2}{\eta_TT\mu^2q}\ln(m+T)\nonumber \\
 & = \bigg( \frac{\rho(F(x_1) - F^*)}{k\gamma} + \frac{9b_1L^2_f\Delta^2_1}{k\lambda\mu} + \frac{2\sigma^2}{qk\mu^2} + \frac{2m\sigma^2}{k\mu^2q}\ln(m+T)\bigg)\frac{(m+T)^{1/2}}{T},
\end{align}
where the second inequality holds by the above inequality \eqref{eq:H6}.
Let $G = \frac{F(x_1) - F^*}{k\gamma\rho} + \frac{9b_1L^2_f\Delta^2_1}{k\lambda\mu\rho^2} + \frac{2\sigma^2}{q k\mu^2\rho^2} + \frac{2m\sigma^2}{q k\mu^2\rho^2}\ln(m+T)$,
we have
\begin{align} \label{eq:66}
 \frac{1}{T} \sum_{t=1}^T \mathbb{E}\big[ \frac{L_f^2}{4\rho^2}\|y^*(x_t)-y_t\|^2 + \frac{1}{4\rho^2}\|\nabla_{x} f(x_t,y_t) -v_t\|^2 +\frac{1}{4\gamma^2}\|\tilde{x}_{t+1}-x_t\|^2 \big]  \leq \frac{G}{T}(m+T)^{1/2}.
\end{align}
According to Jensen's inequality, we have
\begin{align}
 &  \frac{1}{T} \sum_{t=1}^T \mathbb{E}\big[\frac{L_f}{2\rho}\|y^*(x_t)-y_t\| + \frac{1}{2\rho}\|\nabla_{x} f(x_t,y_t) -v_t\| +\frac{1}{2\gamma}\|\tilde{x}_{t+1}-x_t\| \big] \nonumber \\
 & \leq \bigg( \frac{3}{T} \sum_{t=1}^T \mathbb{E}\big[ \frac{L_f^2}{4\rho^2}\|y^*(x_t)-y_t\|^2 + \frac{1}{4\rho^2}\|\nabla_{x} f(x_t,y_t) -v_t\|^2 +\frac{1}{4\gamma^2}\|\tilde{x}_{t+1}-x_t\|^2 \big] \bigg)^{1/2}
 \nonumber \\
 & \leq \frac{\sqrt{3G}}{T^{1/2}}(m+T)^{1/4} \leq \frac{\sqrt{3G}m^{1/4}}{T^{1/2}} + \frac{\sqrt{3G}}{T^{1/4}},
\end{align}
where the last inequality is due to $(a+b)^{1/4} \leq a^{1/4} + b^{1/4}$ for all $a,b>0$. Thus, we have
\begin{align} \label{eq:67}
\frac{1}{T} \sum_{t=1}^T \mathbb{E}\big[  \frac{L_f}{\rho}\|y^*(x_t)-y_t\| + \frac{1}{\rho}\|\nabla_{x} f(x_t,y_t) -v_t\| +\frac{1}{\gamma}\|\tilde{x}_{t+1}-x_t\| \big]  \leq \frac{2\sqrt{3G}m^{1/4}}{T^{1/2}} + \frac{2\sqrt{3G}}{T^{1/4}}.
\end{align}

Let $\phi_t(x)=\frac{1}{2}x^TA_t x$, according to Assumption \ref{ass:4},
$\phi_t(x)$ is $\rho$-strongly convex. Then we define a prox-function (i.e., Bregman distance) associated with $\phi_t(x)$ as in  \cite{censor1981iterative,censor1992proximal,ghadimi2016mini}, defined as
\begin{align}
 D_t(x,x_t) = \phi_t(x) - \big[ \phi_t(x_t) + \langle\nabla \phi_t(x_t), x-x_t\rangle\big] = \frac{1}{2}(x-x_t)^TA_t(x-x_t).
\end{align}
The line 5 of Algorithms \ref{alg:1} is equivalent to the following generalized projection problem
\begin{align}
 \tilde{x}_{t+1} = \arg\min_{x\in \mathcal{X}}\big\{ \langle v_t, x\rangle + \frac{1}{\gamma}D_t(x,x_t)\big\}.
\end{align}
As in \cite{ghadimi2016mini}, we define a generalized projected gradient $\mathcal{G}_{\mathcal{X}}(x_t,v_t,\gamma)=\frac{1}{\gamma}(x_t-\tilde{x}_{t+1})$.
At the same time, we define a gradient mapping $\mathcal{G}_{\mathcal{X}}(x_t,\nabla F(x_t),\gamma)=\frac{1}{\gamma}(x_t-x^*_{t+1})$, where
\begin{align}
 x^*_{t+1} = \arg\min_{x\in \mathcal{X}}\big\{ \langle \nabla F(x_t), x\rangle + \frac{1}{\gamma}D_t(x,x_t)\big\}.
\end{align}
Since $F(x_t) = f(x_t,y^*(x_t)) = \min_{y\in \mathcal{Y}} f(x_t,y)$, by Assumption \ref{ass:5}, we have
\begin{align}
 \|\nabla F(x_t) - v_t\| & = \|\nabla_x f(x_t,y^*(x_t)) - v_t\| \nonumber \\
& = \|\nabla_x f(x_t,y^*(x_t)) - \nabla_x f(x_t,y_t) + \nabla_x f(x_t,y_t)- v_t\| \nonumber \\
&\leq \|\nabla_x f(x_t,y^*(x_t)) - \nabla_x f(x_t,y_t)\| + \|\nabla_x f(x_t,y_t)- v_t\| \nonumber \\
&\leq L_f \|y^*(x_t) - y_t\| + \|\nabla_x f(x_t,y_t)- v_t\|.
\end{align}
According to Proposition 1 in \cite{ghadimi2016mini}, we have $\|\mathcal{G}_{\mathcal{X}}(x_t,\nabla F(x_t),\gamma) - \mathcal{G}_{\mathcal{X}}(x_t,v_t,\gamma) \| \leq  \frac{1}{\rho}\|v_t - \nabla F(x_t)\|$.
Let $ \mathcal{M}_t = \frac{1}{\gamma}\|x_t - \tilde{x}_{t+1}\|
+ \frac{1}{\rho}\big( L_f \|y^*(x_t) - y_t\| + \|\nabla_x f(x_t,y_t)- v_t\|\big)$,
we have
\begin{align}  \label{eq:73}
\|\mathcal{G}_{\mathcal{X}}(x_t,\nabla F(x_t),\gamma)\|
& \leq \|\mathcal{G}_{\mathcal{X}}(x_t,v_t,\gamma)\| + \|\mathcal{G}_{\mathcal{X}}(x_t,\nabla F(x_t),\gamma) - \mathcal{G}_{\mathcal{X}}(x_t,v_t,\gamma) \|  \\
&\leq \|\mathcal{G}_{\mathcal{X}}(x_t,v_t,\gamma)\|
+ \frac{1}{\rho}\|\nabla  F(x_t)-v_t\|  \nonumber \\
& \leq \frac{1}{\gamma}\|x_t - \tilde{x}_{t+1}\|
+ \frac{1}{\rho}\big( L_f \|y^*(x_t) - y_t\| + \|\nabla_x f(x_t,y_t)- v_t\|\big)  = \mathcal{M}_t, \nonumber
\end{align}
where the second inequality holds by the above inequality $\| \nabla F(x_t) - v_t\| \leq L_f \|y^*(x_t) - y_t\| + \|\nabla_x f(x_t,y_t)- v_t\|$.

According to the above inequalities \eqref{eq:73} and \eqref{eq:67}, we have
\begin{align} \label{eq:74}
 \frac{1}{T} \sum_{t=1}^T\mathbb{E}\|\mathcal{G}_{\mathcal{X}}(x_t,\nabla F(x_t),\gamma)\|\leq \frac{1}{T} \sum_{t=1}^T \mathbb{E}\big[\mathcal{M}_t\big]  \leq \frac{2\sqrt{3G}m^{1/4}}{T^{1/2}} + \frac{2\sqrt{3G}}{T^{1/4}}.
\end{align}

\end{proof}

\begin{theorem} \label{th:A2}
(Restatement of Theorem 2)
Assume that the sequence $\{x_t,y_t\}_{t=1}^T$ be generated from the Algorithm \ref{alg:1}. When $\mathcal{X} = \mathbb{R}^{d_1}$, and given $B_t=b_tI_{d_2} \ (\hat{b} \geq b_t \geq b>0)$ for all $t\geq 1$,  $\eta_t=\frac{k}{(m+t)^{1/2}}$ for all $t\geq 0$, $\alpha_{t+1}=c_1\eta_t$, $\beta_{t+1}=c_2\eta_t$, $m\geq \max\big(k^2, (c_1k)^2, (c_2k)^2\big)$, $k>0$, $\frac{9\mu^2}{4} \leq c_1 \leq \frac{m^{1/2}}{k}$, $\frac{75L^2_f}{2} \leq c_2 \leq \frac{m^{1/2}}{k}$, $0< \gamma \leq \min\big(\frac{15\sqrt{2}\lambda\mu^2\rho}{2\sqrt{400L^2_f\lambda^2+24\mu^2\lambda^2+16875\hat{b}^2\kappa^2L^2_f\mu^2}}, \frac{m^{1/2}\rho}{4Lk}\big)$ and $0<\lambda\leq \min\big(\frac{405bL^2_f\mu^{3/2}}{8\sqrt{50L^2_f+9\mu^2}},\frac{b}{6L_f}\big)$, we have
\begin{align}
 \frac{1}{T} \sum_{t=1}^T \mathbb{E}  \|\nabla  F(x_t)\| \leq \frac{\sqrt{\frac{1}{T}\sum_{t=1}^T\mathbb{E}\|A_t\|^2}}{\rho}\bigg(\frac{2\sqrt{3G'}m^{1/4}}{T^{1/2}} + \frac{2\sqrt{3G'}}{T^{1/4}}\bigg),
\end{align}
where $G'=\frac{\rho(F(x_1) - F^*)}{k\gamma} + \frac{9b_1L^2_f\Delta^2_1}{k\lambda\mu} + \frac{2\sigma^2}{q k\mu^2} + \frac{2m\sigma^2}{q k\mu^2}\ln(m+T)$.
\end{theorem}

\begin{proof}
Since $F(x_t) = f(x_t,y^*(x_t)) = \min_{y\in \mathcal{Y}} f(x_t,y)$, we have
\begin{align}
\| \nabla F(x_t) - v_t\| & = \|\nabla_x f(x_t,y^*(x_t)) - v_t\| = \|\nabla_x f(x_t,y^*(x_t)) - \nabla_x f(x_t,y_t)+ \nabla_x f(x_t,y_t)- v_t\| \nonumber \\
& \leq \|\nabla_x f(x_t,y^*(x_t)) - \nabla_x f(x_t,y_t)\| + \|\nabla_x f(x_t,y_t)- v_t\| \nonumber \\
& \leq L_f \|y^*(x_t) - y_t\| + \|\nabla_x f(x_t,y_t)- v_t\|.
\end{align}
Then we have
\begin{align}  \label{eq:11A}
\mathcal{M}_t & = \frac{1}{\gamma}\|x_t - \tilde{x}_{t+1}\|
+ \frac{1}{\rho}\big( L_f \|y^*(x_t) - y_t\|  + \|\nabla_x f(x_t,y_t)- v_t\|\big) \nonumber \\
& \geq \frac{1}{\gamma}\|x_t - \tilde{x}_{t+1}\|
+ \frac{1}{\rho}\|\nabla  F(x_t)-v_t\| \nonumber \\
& \mathop{=}^{(i)} \|A_t^{-1}v_t\|
+ \frac{1}{\rho}\|\nabla  F(x_t)-v_t\| \nonumber \\
& = \frac{1}{\|A_t\|}\|A_t\|\|A_t^{-1}v_t\|
+ \frac{1}{\rho}\|\nabla  F(x_t)-v_t\| \nonumber \\
& \geq \frac{1}{\|A_t\|}\|v_t\|
+ \frac{1}{\rho}\|\nabla  F(x_t)-v_t\| \nonumber \\
& \mathop{\geq}^{(ii)} \frac{1}{\|A_t\|}\|v_t\|
+ \frac{1}{\|A_t\|}\|\nabla  F(x_t)-v_t\| \nonumber \\
& \geq \frac{1}{\|A_t\|}\|\nabla  F(x_t)\|
\end{align}
where the equality $(i)$ holds by $\tilde{x}_{t+1} = x_t - \gamma A^{-1}_tv_t$ that can be easily obtained from the line 5 of Algorithm \ref{alg:1} when $\mathcal{X}=\mathbb{R}^{d_1}$, and the inequality $(ii)$ holds by $\|A_t\| \geq \rho$ for all $t\geq1$ due to Assumption \ref{ass:4}.
Then we have
 \begin{align}
  \|\nabla F(x_t)\| \leq \mathcal{M}_t \|A_t\|.
 \end{align}
By using Cauchy-Schwarz inequality, we have
 \begin{align} \label{eq:78}
  \frac{1}{T}\sum_{t=1}^T\mathbb{E}\|\nabla F(x_t)\| \leq \frac{1}{T}\sum_{t=1}^T\mathbb{E}\big[\mathcal{M}_t \|A_t\|\big] \leq \sqrt{\frac{1}{T}\sum_{t=1}^T\mathbb{E}[\mathcal{M}_t^2]} \sqrt{\frac{1}{T}\sum_{t=1}^T\mathbb{E}\|A_t\|^2}.
 \end{align}
According to the above inequality \eqref{eq:66} and $ \mathcal{M}_t = \frac{1}{\gamma}\|x_t - \tilde{x}_{t+1}\|
+ \frac{1}{\rho}\big( L_f \|y^*(x_t) - y_t\| + \|\nabla_x f(x_t,y_t)- v_t\|\big)$,
we have
\begin{align} \label{eq:79}
 \frac{1}{T}\sum_{t=1}^T\mathbb{E}\big[\mathcal{M}_t^2\big] & \leq \frac{1}{T}\sum_{t=1}^T\big[\frac{3L_f^2}{\rho^2}\|y^*(x_t)-y_t\|^2 + \frac{3}{\rho^2}\|\nabla_{x} f(x_t,y_t) -v_t\|^2 +\frac{3}{\gamma^2}\|\tilde{x}_{t+1}-x_t\|^2 \big]  \nonumber \\
 & \leq \frac{12G}{T}(m+T)^{1/2}.
\end{align}
By combining the above inequalities \eqref{eq:78} and \eqref{eq:79}, we have
\begin{align}
\frac{1}{T}\sum_{t=1}^T\mathbb{E}\|\nabla F(x_t)\| \leq  \sqrt{\frac{1}{T}\sum_{t=1}^T\mathbb{E}\|A_t\|^2}\frac{2\sqrt{3G}}{T^{1/2}}(m+T)^{1/4}.
\end{align}
Let $G' = \rho^2 G =\frac{\rho(F(x_1) - F^*)}{k\gamma} + \frac{9b_1L^2_f\Delta^2_1}{k\lambda\mu} + \frac{2\sigma^2}{q k\mu^2} + \frac{2m\sigma^2}{q k\mu^2}\ln(m+T)$,
we have
\begin{align}
\frac{1}{T} \sum_{t=1}^T \mathbb{E}\|\nabla  F(x_t)\| \leq \frac{\sqrt{\frac{1}{T}\sum_{t=1}^T\mathbb{E}\|A_t\|^2}}{\rho}\bigg(\frac{2\sqrt{3G'}m^{1/4}}{T^{1/2}} + \frac{2\sqrt{3G'}}{T^{1/4}}\bigg).
\end{align}

\end{proof}

\subsection{ Convergence Analysis of the VR-AdaGDA Algorithm }
\label{Appendix:A2}
In the subsection, we study the convergence properties of the VR-AdaGDA algorithm for solving the minimax problem \eqref{eq:1}.
We first provide a useful lemma.

\begin{lemma} \label{lem:F1}
Suppose the stochastic gradients $v_t$ and $w_t$ be generated from Algorithm \ref{alg:2}, we have
\begin{align} \label{eq:F1}
 \mathbb{E} \|\nabla_{x} f(x_{t+1},y_{t+1}) - v_{t+1}\|^2 & \leq (1-\alpha_{t+1})\mathbb{E} \|\nabla_{x} f(x_t,y_t) -v_t\|^2 + \frac{2\alpha_{t+1}^2\sigma^2}{q} \nonumber \\
 & \quad + \frac{4L^2_f\eta^2_t}{q}\mathbb{E}\big(\|\tilde{x}_{t+1}-x_t\|^2 + \|\tilde{y}_{t+1}-y_t\|^2\big),
\end{align}
\begin{align} \label{eq:F2}
\mathbb{E} \|\nabla_{y} f(x_{t+1},y_{t+1}) - w_{t+1}\|^2 & \leq (1-\beta_{t+1})\mathbb{E} \|\nabla_{y} f(x_t,y_t) -w_t\|^2 + \frac{2\beta_{t+1}^2\sigma^2}{q} \nonumber \\
 & \quad + \frac{4L^2_f\eta^2_t}{q}\mathbb{E}\big(\|\tilde{x}_{t+1} - x_t\|^2 + \|\tilde{y}_{t+1}-y_t\|^2\big).
\end{align}
\end{lemma}

\begin{proof}
We first prove the inequality \eqref{eq:F1}.
According to the definition of $v_t$ in Algorithm \ref{alg:2}, we have
\begin{align}
v_{t+1} - v_t & = -\alpha_{t+1}v_t + (1-\alpha_{t+1})\big(\nabla_{x} f(x_{t+1},y_{t+1};\mathcal{B}_{t+1}) - \nabla_{x} f(x_t,y_t;\mathcal{B}_{t+1})\big) \nonumber \\
& \quad + \alpha_{t+1}\nabla_{x} f(x_{t+1},y_{t+1};\mathcal{B}_{t+1}).
\end{align}
Then we have
\begin{align} \label{eq:F3}
& \mathbb{E} \|\nabla_{x} f(x_{t+1},y_{t+1}) - v_{t+1}\|^2 \nonumber \\
&= \mathbb{E} \|\nabla_{x} f(x_{t+1},y_{t+1}) - v_t - (v_{t+1}-v_t)\|^2 \\
& = \mathbb{E} \|\nabla_{x} f(x_{t+1},y_{t+1}) - v_t + \alpha_{t+1}v_t- \alpha_{t+1}\nabla_{x} f(x_{t+1},y_{t+1};\mathcal{B}_{t+1})
- (1-\alpha_{t+1})\big(\nabla_{x} f(x_{t+1},y_{t+1};\mathcal{B}_{t+1}) \nonumber \\
&\quad - \nabla_{x} f(x_t,y_t;\mathcal{B}_{t+1})\big) \|^2 \nonumber \\
& = \mathbb{E} \|(1-\alpha_{t+1})(\nabla_{x} f(x_t,y_t) -v_t) + (1-\alpha_{t+1})\big(\nabla_{x} f(x_{t+1},y_{t+1})
-\nabla_{x} f(x_t,y_t)-\nabla_{x} f(x_{t+1},y_{t+1};\mathcal{B}_{t+1}) \nonumber \\
& \quad + \nabla_{x} f(x_t,y_t;\mathcal{B}_{t+1})\big)+ \alpha_{t+1}\big(\nabla_{x} f(x_{t+1},y_{t+1})- \nabla_{x} f(x_{t+1},y_{t+1};\mathcal{B}_{t+1})\big)\|^2 \nonumber \\
& = (1-\alpha_{t+1})^2\mathbb{E}\|\nabla_{x} f(x_t,y_t) - v_t\|^2 + \alpha_{t+1}^2 \mathbb{E} \|\nabla_{x} f(x_{t+1},y_{t+1})- \nabla_{x} f(x_{t+1},y_{t+1};\mathcal{B}_{t+1})\|^2\nonumber \\
& \quad +(1- \alpha_{t+1})^2\mathbb{E} \|\nabla_{x} f(x_{t+1},y_{t+1})
- \nabla_{x} f(x_t,y_t)  -  \nabla_{x} f(x_{t+1},y_{t+1};\mathcal{B}_{t+1}) + \nabla_{x} f(x_t,y_t;\mathcal{B}_{t+1})\|^2
\nonumber \\
& \quad+ 2\alpha_{t+1}(1-\alpha_{t+1})\big\langle\nabla_{x} f(x_{t+1},y_{t+1})
- \nabla_{x} f(x_t,y_t)  -  \nabla_{x} f(x_{t+1},y_{t+1};\mathcal{B}_{t+1}) + \nabla_{x} f(x_t,y_t;\mathcal{B}_{t+1}), \nonumber \\
& \quad \quad \nabla_{x} f(x_{t+1},y_{t+1})  -  \nabla_{x} f(x_{t+1},y_{t+1};\mathcal{B}_{t+1})\big\rangle \nonumber \\
& \leq (1-\alpha_{t+1})^2\mathbb{E} \|\nabla_{x} f(x_t,y_t) -v_t\|^2 + 2\alpha_{t+1}^2 \mathbb{E} \|\nabla_{x} f(x_{t+1},y_{t+1})- \nabla_{x} f(x_{t+1},y_{t+1};\mathcal{B}_{t+1})\|^2  \nonumber \\
& \quad + 2(1-\alpha_{t+1})^2\mathbb{E} \|\nabla_{x} f(x_{t+1},y_{t+1})  -  \nabla_{x} f(x_t,y_t) -  \nabla_{x} f(x_{t+1},y_{t+1};\mathcal{B}_{t+1})  +  \nabla_{x} f(x_t,y_t;\mathcal{B}_{t+1})\|^2 \nonumber \\
& \leq (1-\alpha_{t+1})^2 \mathbb{E} \|\nabla_{x} f(x_t,y_t) -v_t\|^2 + \frac{2\alpha_{t+1}^2\sigma^2}{q} \nonumber \\
& \quad + \frac{2(1-\alpha_{t+1})^2}{q}\underbrace{ \mathbb{E} \| \nabla_{x} f(x_{t+1},y_{t+1};\xi_{t+1}) -
\nabla_{x} f(x_t,y_t;\xi_{t+1})\|^2}_{=T_1}, \nonumber
\end{align}
where the fourth equality follows by $\mathbb{E}_{\mathcal{B}_{t+1}}[\nabla_{x} f(x_{t+1},y_{t+1};\mathcal{B}_{t+1})]=\nabla_{x} f(x_{t+1},y_{t+1})$ and $ \mathbb{E}_{\mathcal{B}_{t+1}}[\nabla_{x} f(x_{t+1},y_{t+1};\mathcal{B}_{t+1}) - \nabla_{x} f(x_t,y_t;\mathcal{B}_{t+1})]=\nabla_{x} f(x_{t+1},y_{t+1}) - \nabla_{x} f(x_t,y_t)$; the first inequality holds by Young's inequality; the last inequality is due to Lemma \ref{lem:A4} and Assumption \ref{ass:1}.

According to Assumption \ref{ass:6}, we have
\begin{align} \label{eq:F4}
T_1 & = \mathbb{E} \big\| \nabla_{x} f(x_{t+1},y_{t+1};\xi_{t+1}) - \nabla_{x} f(x_t,y_t;\xi_{t+1})\big\|^2  \nonumber \\
& = \mathbb{E} \big\| \nabla_{x} f(x_{t+1},y_{t+1};\xi_{t+1}) - \nabla_{x} f(x_t,y_{t+1};\xi_{t+1}) + \nabla_{x} f(x_t,y_{t+1};\xi_{t+1})- \nabla_{x} f(x_t,y_t;\xi_{t+1})\big\|^2  \nonumber \\
& \leq 2L_f^2 \mathbb{E}\|x_{t+1}-x_t\|^2 + 2L_f^2\mathbb{E}\|y_{t+1}-y_t\|^2  \nonumber \\
& = 2L_f^2 \eta_t^2\mathbb{E}\|\tilde{x}_{t+1}-x_t\|^2 + 2L_f^2\eta_t^2\mathbb{E}\|\tilde{y}_{t+1}-y_t\|^2.
\end{align}
Plugging the above inequality \eqref{eq:F4} into \eqref{eq:F3}, we obtain
\begin{align}
 \mathbb{E} \|\nabla_{x} f(x_{t+1},y_{t+1}) - v_{t+1}\|^2 & \leq (1-\alpha_{t+1})^2 \mathbb{E} \|\nabla_{x} f(x_t,y_t) -v_t\|^2 + \frac{2\alpha_{t+1}^2\sigma^2}{q} \nonumber \\
 & \quad + \frac{4(1-\alpha_{t+1})^2L^2_f\eta^2_t}{q}\mathbb{E}\big(\|\tilde{x}_{t+1}-x_t\|^2 + \|\tilde{y}_{t+1}-y_t\|^2\big) \nonumber \\
 & \leq (1-\alpha_{t+1})\mathbb{E} \|\nabla_{x} f(x_t,y_t) -v_t\|^2 + \frac{2\alpha_{t+1}^2\sigma^2}{q} \nonumber \\
 & \quad + \frac{4L^2_f\eta^2_t}{q}\mathbb{E}\big(\|\tilde{x}_{t+1}-x_t\|^2 + \|\tilde{y}_{t+1}-y_t\|^2\big),
\end{align}
where the last inequality holds by $0<\alpha_{t+1}\leq 1$.

Similarly, we have
\begin{align}
 \mathbb{E} \|\nabla_{y} f(x_{t+1},y_{t+1}) - w_{t+1}\|^2 & \leq (1-\beta_{t+1}) \mathbb{E} \|\nabla_{y} f(x_t,y_t) -w_t\|^2 + \frac{2\beta_{t+1}^2\sigma^2}{q} \nonumber \\
 & \quad + \frac{4L^2_f\eta^2_t}{q}\mathbb{E}\big(\|\tilde{x}_{t+1} - x_t\|^2 + \|\tilde{y}_{t+1}-y_t\|^2\big).
\end{align}

\end{proof}

\begin{theorem} \label{th:A3}
(Restatement of Theorem 3)
Suppose the sequence $\{x_t,y_t\}_{t=1}^T$ be generated from Algorithm \ref{alg:2}. When $\mathcal{X}\subset \mathbb{R}^{d_1}$, and given $B_t=b_tI_{d_2} \ (\hat{b}\geq b_t\geq b >0)$, $\eta_t=\frac{k}{(m+t)^{1/3}}$, $\alpha_{t+1}=c_1\eta^2_t$, $\beta_{t+1}=c_2\eta^2_t$, $c_1 \geq \frac{2}{3k^3} + \frac{9\mu^2}{4}$ and $c_2 \geq \frac{2}{3k^3} + \frac{75L^2_f}{2}$, $m \geq \max\big( k^3, (c_1k)^3, (c_2k)^3\big)$, $0< \lambda \leq \min\big(\frac{27\mu bq}{32},\frac{b}{6L_f}\big)$ and
$0< \gamma \leq \min\big( \frac{\rho\lambda\mu\sqrt{q}}{L_f\sqrt{32\lambda^2+150q\kappa^2\hat{b}^2}}, \frac{m^{1/3}\rho}{2Lk} \big)$,
we have
\begin{align}
 \frac{1}{T}\sum_{t=1}^T\mathbb{E}\|\mathcal{G}_{\mathcal{X}}(x_t,\nabla F(x_t),\gamma)\| \leq \frac{2\sqrt{3M}m^{1/6}}{T^{1/2}} + \frac{2\sqrt{3M}}{T^{1/3}},
\end{align}
where $M = \frac{F(x_1) - F^*}{T\gamma k\rho} + \frac{9L^2_f b_1}{k\lambda\mu\rho^2}\Delta_1^2  + \frac{2\sigma^2m^{1/3}}{k^2q\mu^2\rho^2} + \frac{2k^2(c_1^2+c_2^2)\sigma^2}{q\mu^2\rho^2}\ln(m+T)$ and $\Delta_1^2 = \|y_1 - y^*(x_1)\|^2$.
\end{theorem}

\begin{proof}
Since $\eta_t$ is decreasing and $m\geq k^3$, we have $\eta_t \leq \eta_0 = \frac{k}{m^{1/3}} \leq 1$ and $\gamma \leq \frac{\rho}{2L\eta_0} =\frac{m^{1/3}\rho}{2Lk} \leq \frac{1}{2L\eta_t}$ for any $t\geq 0$.
Due to $0 < \eta_t \leq 1$ and $m\geq \max\big( (c_1k)^3, (c_2k)^3 \big)$, we have $\alpha_t = c_1\eta_t^2 \leq c_1\eta_t \leq \frac{c_1k}{m^{1/3}}\leq 1$ and $\beta_t = c_2\eta_t^2 \leq c_2\eta_t \leq \frac{c_2k}{m^{1/3}}\leq 1$.
Then we consider the upper bound of the following term:
\begin{align}
& \frac{1}{\eta_t}\mathbb{E} \|\nabla_{x} f(x_{t+1},y_{t+1}) - v_{t+1}\|^2 - \frac{1}{\eta_{t-1}}\mathbb{E} \|\nabla_{x} f(x_t,y_t) - v_t\|^2  \\
& \leq \big(\frac{1 - \alpha_{t+1}}{\eta_t}  -  \frac{1}{\eta_{t-1}}\big)\mathbb{E} \|\nabla_{x} f(x_t,y_t) -v_t\|^2  +  \frac{4L^2_f\eta_{t}}{q}\mathbb{E}\big(\|\tilde{x}_{t+1}-x_t\|^2 + \|\tilde{y}_{t+1}-y_t\|^2\big)  +  \frac{2\alpha_{t+1}^2\sigma^2}{q\eta_t}\nonumber \\
& = \big(\frac{1}{\eta_t}  -  \frac{1}{\eta_{t-1}} - c_1\eta_t\big)\mathbb{E} \|\nabla_{x} f(x_t,y_t) -v_t\|^2 + \frac{4L^2_f\eta_{t}}{q}\mathbb{E}\big(\|\tilde{x}_{t+1}-x_t\|^2 + \|\tilde{y}_{t+1}-y_t\|^2\big)  +  \frac{2c^2_1\eta_t^3\sigma^2}{q}, \nonumber
\end{align}
where the second inequality is due to $0<\alpha_{t+1}\leq 1$.
Similarly, we have
\begin{align}
& \frac{1}{\eta_t}\mathbb{E} \|\nabla_{y} f(x_{t+1},y_{t+1}) - w_{t+1}\|^2
-  \frac{1}{\eta_{t-1}}\mathbb{E} \|\nabla_{y} f(x_t,y_t) - w_t\|^2  \\
& \leq \big(\frac{1}{\eta_t} - \frac{1}{\eta_{t-1}} - c_2\eta_t\big)\mathbb{E} \|\nabla_{y} f(x_t,y_t) -w_t\|^2 + \frac{4L^2_f\eta_{t}}{q}\mathbb{E}\big(\|\tilde{x}_{t+1}-x_t\|^2 + \|\tilde{y}_{t+1}-y_t\|^2\big) +
 \frac{2c^2_2\eta_t^3\sigma^2}{q}. \nonumber
\end{align}
By $\eta_t = \frac{k}{(m+t)^{1/3}}$, we have
\begin{align}
\frac{1}{\eta_t} - \frac{1}{\eta_{t-1}} &= \frac{1}{k}\big( (m+t)^{\frac{1}{3}} - (m+t-1)^{\frac{1}{3}}\big) \nonumber \\
& \leq \frac{1}{3k(m+t-1)^{2/3}} = \frac{2^{2/3}}{3k\big(2(m+t-1)\big)^{2/3}} \nonumber \\
& \leq \frac{2^{2/3}}{3k(m+t)^{2/3}} = \frac{2^{2/3}}{3k^3}\frac{k^2}{(m+t)^{2/3}}= \frac{2^{2/3}}{3k^3}\eta_t^2 \leq \frac{2}{3k^3}\eta_t,
\end{align}
where the first inequality holds by the concavity of function $f(x)=x^{1/3}$, \emph{i.e.}, $(x+y)^{1/3}\leq x^{1/3} + \frac{y}{3x^{2/3}}$, and
the last inequality is due to $0<\eta_t\leq 1$.

Let $c_1 \geq \frac{2}{3k^3} + \frac{9\mu^2}{4}$, we have
\begin{align} \label{eq:G1}
& \frac{1}{\eta_t}\mathbb{E} \|\nabla_{x} f(x_{t+1},y_{t+1}) - v_{t+1}\|^2 - \frac{1}{\eta_{t-1}}\mathbb{E} \|\nabla_{x} f(x_t,y_t) - v_t\|^2  \\
& \leq -\frac{9\mu^2\eta_t}{4}\mathbb{E} \|\nabla_{x} f(x_t,y_t) -v_t\|^2 + \frac{4L^2_f\eta_{t}}{q}\mathbb{E}\big(\|\tilde{x}_{t+1}-x_t\|^2 + \|\tilde{y}_{t+1}-y_t\|^2\big) + \frac{2c^2_1\eta_t^3\sigma^2}{q}. \nonumber
\end{align}
Let $c_2 \geq \frac{2}{3k^3} + \frac{75L^2_f}{2}$, we have
\begin{align} \label{eq:G2}
& \frac{1}{\eta_t}\mathbb{E} \|\nabla_{y} f(x_{t+1},y_{t+1}) - w_{t+1}\|^2
-  \frac{1}{\eta_{t-1}}\mathbb{E} \|\nabla_{y} f(x_t,y_t) - w_t\|^2  \\
& \leq - \frac{75L^2_f\eta_t}{2}\mathbb{E} \|\nabla_{y} f(x_t,y_t) -w_t\|^2 + \frac{4L^2_f\eta_{t}}{q}\mathbb{E}\big(\|\tilde{x}_{t+1}-x_t\|^2 + \|\tilde{y}_{t+1}-y_t\|^2\big) + \frac{2c^2_2\eta_t^3\sigma^2}{q}. \nonumber
 \end{align}

According to Lemma \ref{lem:D1}, we have
\begin{align} \label{eq:G3}
F(x_{t+1}) - F(x_t) \leq \frac{2\gamma L_f^2\eta_t}{\rho}\|y^*(x_t)-y_t\|^2 + \frac{2\gamma\eta_t}{\rho}\|\nabla_{x} f(x_t,y_t) -v_t \|^2 -\frac{\rho\eta_t}{2\gamma}\|\tilde{x}_{t+1}-x_t\|^2.
\end{align}
According to Lemma \ref{lem:E1}, we have
\begin{align} \label{eq:G4}
\|y_{t+1} - y^*(x_{t+1})\|^2 - \|y_t -y^*(x_t)\|^2  & \leq -\frac{\eta_t\mu\lambda}{4b_t} \|y_t -y^*(x_t)\|^2 -\frac{3\eta_t}{4} \|\tilde{y}_{t+1}-y_t\|^2 \nonumber \\
& \quad + \frac{25\eta_t\lambda}{6\mu b_t} \|\nabla_y f(x_t,y_t)-w_t\|^2 + \frac{25\kappa^2\eta_tb_t}{6\mu\lambda}\|\tilde{x}_{t+1} - x_t\|^2.
\end{align}

Next, we define a Lyapunov function, for any $t\geq 1$
\begin{align}
\Phi_t = \mathbb{E}\big[ F(x_t) + \frac{9\gamma L^2_fb_t}{\rho\lambda\mu}\|y_t - y^*(x_t)\|^2 + \frac{\gamma}{\rho\mu^2}\big(\frac{1}{\eta_{t-1}}\|\nabla_{x} f(x_t,y_t)-v_t\|^2 + \frac{1}{\eta_{t-1}}\|\nabla_{y} f(x_t,y_t)-w_t\|^2 \big) \big].
\end{align}
Then we have
\begin{align}
& \Phi_{t+1}- \Phi_t \nonumber \\
& = \mathbb{E}\big[F(x_{t+1})  -  F(x_t)\big] \!+\! \frac{9\gamma L^2_fb_t}{\rho\lambda\mu}
\big(\mathbb{E}\|y_{t+1}-y^*(x_{t+1})\|^2 \!-\! \mathbb{E} \|y_t-y^*(x_t)\|^2 \big)
 \!+\! \frac{\gamma}{\rho\mu^2} \big(\frac{1}{\eta_t}\mathbb{E}\|\nabla_{x} f(x_{t+1},y_{t+1})-v_{t+1}\|^2 \nonumber \\
& \quad  -  \frac{1}{\eta_{t-1}}\mathbb{E}\|\nabla_{x} f(x_t,y_t)-v_t\|^2  +  \frac{1}{\eta_t}\mathbb{E}\|\nabla_{y} f(x_{t+1},y_{t+1})-w_{t+1}\|^2
 -  \frac{1}{\eta_{t-1}}\mathbb{E}\|\nabla_{y} f(x_t,y_t)-w_t\|^2 \big) \nonumber \\
& \leq \frac{2\gamma L_f^2\eta_t}{\rho}\mathbb{E}\|y^*(x_t)-y_t\|^2 + \frac{2\gamma\eta_t}{\rho}\mathbb{E}\|\nabla_{x} f(x_t,y_t) -v_t \|^2 -\frac{\rho\eta_t}{2\gamma}\mathbb{E}\|\tilde{x}_{t+1}-x_t\|^2 \nonumber \\
& \quad  + \frac{9\gamma L^2_fb_t}{\rho\lambda\mu} \bigg( \!-\!\frac{\eta_t\mu\lambda}{4b_t}\mathbb{E}\|y_t -y^*(x_t)\|^2 \!-\!\frac{3\eta_t}{4} \mathbb{E}\|\tilde{y}_{t+1}-y_t\|^2 \!+\! \frac{25\eta_t\lambda}{6\mu b_t}\mathbb{E} \|\nabla_y f(x_t,y_t)-w_t\|^2 \!+\! \frac{25\kappa^2\eta_tb_t}{6\mu\lambda}\mathbb{E}\|\tilde{x}_{t+1} - x_t\|^2 \bigg)  \nonumber \\
& \quad -\frac{9\gamma\eta_t}{4\rho} \mathbb{E} \|\nabla_x f(x_t,y_t) -v_t\|^2 + \frac{4\gamma L^2_f\eta_{t}}{\rho \mu^2q}\mathbb{E}\big(\|\tilde{x}_{t+1}-x_t\|^2 + \|\tilde{y}_{t+1}-y_t\|^2\big) + \frac{2\gamma c^2_1\eta_t^3\sigma^2}{\rho \mu^2q} \nonumber \\
& \quad -\frac{75L^2_f\gamma}{2\mu^2\rho} \eta_t\mathbb{E}\|\nabla_{y} f(x_t,y_t) -w_t\|^2  + \frac{4\gamma L^2_f\eta_{t}}{\rho \mu^2q}\mathbb{E}\big(\|\tilde{x}_{t+1}-x_t\|^2 + \|\tilde{y}_{t+1}-y_t\|^2\big) + \frac{2\gamma c^2_2\eta_t^3\sigma^2}{\rho \mu^2q} \nonumber \\
& \leq -\frac{\gamma L_f^2\eta_t}{4\rho}\mathbb{E}\|y^*(x_t)-y_t\|^2 - \frac{\gamma\eta_t}{4\rho}\mathbb{E}\|\nabla_{x} f(x_t,y_t) -v_t\|^2 + \frac{2\gamma c^2_1\eta_t^3\sigma^2}{\rho \mu^2q} + \frac{2\gamma c^2_2\eta_t^3\sigma^2}{\rho \mu^2q} \nonumber \\
& \quad -\big(\frac{27b_t\gamma L^2_f}{4\rho\lambda\mu} - \frac{8\gamma L^2_f}{\rho\mu^2q}\big)\eta_t\mathbb{E}\|\tilde{y}_{t+1} - y_{t+1}\|^2 - \big( \frac{\rho}{2\gamma} - \frac{8\gamma L^2_f}{\rho\mu^2q} - \frac{75\gamma L^2_f\kappa^2b_t^2}{2\rho\lambda^2\mu^2}\big)\eta_t\mathbb{E}\|\tilde{x}_{t+1}-x_t\|^2 \nonumber \\
& \leq -\frac{\gamma L_f^2\eta_t}{4\rho}\mathbb{E}\|y^*(x_t)-y_t\|^2  -  \frac{\gamma\eta_t}{4\rho}\mathbb{E}\|\nabla_{x} f(x_t,y_t) -v_t\|^2 - \frac{\rho\eta_t}{4\gamma}\mathbb{E}\|\tilde{x}_{t+1}-x_t\|^2 + \frac{2\gamma c^2_1\eta_t^3\sigma^2}{\rho \mu^2q} + \frac{2\gamma c^2_2\eta_t^3\sigma^2}{\rho \mu^2q},
\end{align}
where the first inequality holds by the above inequalities \eqref{eq:G1}, \eqref{eq:G2}, \eqref{eq:G3} and \eqref{eq:G4};
the last inequality is due to $0< \lambda \leq \frac{27\mu bq}{32} \leq \frac{27\mu b_tq}{32}$ and $0< \gamma \leq \frac{\rho\lambda\mu\sqrt{q}}{L_f\sqrt{32\lambda^2+150q\kappa^2\hat{b}^2}} \leq \frac{\rho\lambda\mu\sqrt{q}}{L_f\sqrt{32\lambda^2+150q\kappa^2b_t^2}}$
for all $t\geq 1$.
Thus, we have
\begin{align} \label{eq:G5}
 & \frac{L_f^2\eta_t}{4}\mathbb{E}\|y^*(x_t)-y_t\|^2 + \frac{\eta_t}{4}\mathbb{E}\|\nabla_{x} f(x_t,y_t) -v_t\|^2 + \frac{\rho^2\eta_t}{4\gamma^2}\mathbb{E}\|\tilde{x}_{t+1}-x_t\|^2 \nonumber \\
 & \leq \frac{\rho(\Phi_t - \Phi_{t+1})}{\gamma} + \frac{2c^2_1\eta_t^3\sigma^2}{\mu^2q} + \frac{2c^2_2\eta_t^3\sigma^2}{\mu^2q}.
\end{align}

Taking average over $t=1,2,\cdots,T$ on both sides of \eqref{eq:G5}, we have
\begin{align}
 & \frac{1}{T} \sum_{t=1}^T \big( \frac{L_f^2\eta_t}{4}\mathbb{E}\|y^*(x_t)-y_t\|^2 + \frac{\eta_t}{4}\mathbb{E}\|\nabla_{x} f(x_t,y_t) -v_t\|^2 +\frac{\rho^2\eta_t}{4\gamma^2}\mathbb{E}\|\tilde{x}_{t+1}-x_t\|^2  \big) \nonumber \\
 & \leq  \sum_{t=1}^T \frac{\rho(\Phi_t - \Phi_{t+1})}{T\gamma} + \frac{1}{T}\sum_{t=1}^T\big(\frac{2c^2_1\eta_t^3\sigma^2}{\mu^2q} + \frac{2c^2_2\eta_t^3\sigma^2}{\mu^2q}\big). \nonumber
\end{align}
Given $x_1 \in \mathcal{X}$, $y_1 \in \mathcal{Y}$ and
 $\Delta_1^2 = \|y_1 - y^*(x_1)\|^2$, we have
\begin{align} \label{eq:G6}
 \Phi_1 &= F(x_1) + \frac{9\gamma L^2_fb_1}{\rho\lambda\mu}\|y_1 - y^*(x_1)\|^2 + \frac{\gamma}{\rho\mu^2\eta_0}\mathbb{E}\|\nabla_{x} f(x_1,y_1)-v_1\|^2 + \frac{\gamma}{\rho\mu^2\eta_{0}}\mathbb{E}\|\nabla_{y} f(x_1,y_1)-w_1\|^2 \nonumber \\
 & \leq F(x_1) + \frac{9\gamma L^2_fb_1}{\rho\lambda\mu}\Delta_1^2  + \frac{2\gamma\sigma^2}{q\rho\mu^2\eta_0},
\end{align}
where the last inequality holds by Assumption \ref{ass:1}.

Since $\eta_t$ is decreasing, i.e., $\eta_T^{-1} \geq \eta_t^{-1}$ for any $0\leq t\leq T$, we have
 \begin{align}
 & \frac{1}{T} \sum_{t=1}^T \mathbb{E}\big(  \frac{L_f^2}{4}\|y^*(x_t)-y_t\|^2 + \frac{1}{4}\|\nabla_{x} f(x_t,y_t) -v_t\|^2 +\frac{\rho^2}{4\gamma^2}\|\tilde{x}_{t+1}-x_t\|^2 \big) \nonumber \\
 & \leq  \sum_{t=1}^T \frac{\rho(\Phi_t - \Phi_{t+1})}{\eta_TT\gamma} + \frac{1}{\eta_TT}\sum_{t=1}^T\big(\frac{2c^2_1\eta_t^3\sigma^2}{\mu^2q} + \frac{2c^2_2\eta_t^3\sigma^2}{\mu^2q}\big) \nonumber \\
 & =  \frac{\rho(\Phi_1 - \Phi_{T+1})}{\eta_TT\gamma} + \frac{2(c_1^2+c_2^2)\sigma^2}{\eta_TTq\mu^2}\sum_{t=1}^T\eta_t^3 \nonumber \\
 & \leq \frac{\rho(F(x_1) - F^*)}{T\eta_T\gamma} + \frac{9L^2_fb_1}{\eta_TT\lambda\mu}\Delta_1^2  + \frac{2\sigma^2}{\eta_TTq\mu^2\eta_0} + \frac{2(c_1^2+c_2^2)\sigma^2}{\eta_TTq\mu^2}\int^T_1\frac{k^3}{m+t}dt \nonumber \\
 & \leq \frac{\rho(F(x_1) - F^*)}{T\eta_T\gamma} + \frac{9L^2_fb_1}{\eta_TT\lambda\mu}\Delta_1^2  + \frac{2\sigma^2}{\eta_TTq\mu^2\eta_0} + \frac{2k^3(c_1^2+c_2^2)\sigma^2}{\eta_TTq\mu^2}\ln(m+T) \nonumber \\
 & = \bigg( \frac{\rho(F(x_1) - F^*)}{T\gamma k} + \frac{9L^2_fb_1}{k\lambda\mu}\Delta_1^2  + \frac{2\sigma^2m^{1/3}}{k^2q\mu^2} + \frac{2k^2(c_1^2+c_2^2)\sigma^2}{q\mu^2}\ln(m+T)\bigg) \frac{(m+T)^{1/3}}{T},
\end{align}
where the second inequality holds by the above inequality \eqref{eq:G6}. Let $M = \frac{F(x_1) - F^*}{T\gamma k\rho} + \frac{9L^2_f b_1}{k\lambda\mu\rho^2}\Delta_1^2  + \frac{2\sigma^2m^{1/3}}{k^2q\mu^2\rho^2} + \frac{2k^2(c_1^2+c_2^2)\sigma^2}{q\mu^2\rho^2}\ln(m+T)$,
we have
\begin{align} \label{eq:102}
 \frac{1}{T} \sum_{t=1}^T \mathbb{E}\big[ \frac{L_f^2}{4\rho^2}\|y^*(x_t)-y_t\|^2 + \frac{1}{4\rho^2}\|\nabla_{x} f(x_t,y_t) -v_t\|^2 +\frac{1}{4\gamma^2}\|\tilde{x}_{t+1}-x_t\|^2 \big]  \leq \frac{M}{T}(m+T)^{1/3}.
\end{align}
According to Jensen's inequality, we have
\begin{align}
 &  \frac{1}{T} \sum_{t=1}^T \mathbb{E}\big[\frac{L_f}{2\rho}\|y^*(x_t)-y_t\| + \frac{1}{2\rho}\|\nabla_{x} f(x_t,y_t) -v_t\| +\frac{1}{2\gamma}\|\tilde{x}_{t+1}-x_t\| \big] \nonumber \\
 & \leq \bigg( \frac{3}{T} \sum_{t=1}^T \mathbb{E}\big[ \frac{L_f^2}{4\rho^2}\|y^*(x_t)-y_t\|^2 + \frac{1}{4\rho^2}\|\nabla_{x} f(x_t,y_t) -v_t\|^2 +\frac{1}{4\gamma^2}\|\tilde{x}_{t+1}-x_t\|^2 \big]\bigg)^{1/2} \nonumber \\
 & \leq \frac{\sqrt{3M}}{T^{1/2}}(m+T)^{1/6} \leq \frac{\sqrt{3M}m^{1/6}}{T^{1/2}} + \frac{\sqrt{3M}}{T^{1/3}},
\end{align}
where the last inequality is due to $(a+b)^{1/6} \leq a^{1/6} + b^{1/6}$. Thus, we have
\begin{align} \label{eq:91}
\frac{1}{T} \sum_{t=1}^T \mathbb{E}\big[ \frac{L_f}{\rho}\|y^*(x_t)-y_t\| + \frac{1}{\rho}\|\nabla_{x} f(x_t,y_t) -v_t\| +\frac{1}{\gamma}\|\tilde{x}_{t+1}-x_t\| \big]  \leq \frac{2\sqrt{3M}m^{1/6}}{T^{1/2}} + \frac{2\sqrt{3M}}{T^{1/3}}.
\end{align}
According to the above inequalities \eqref{eq:73} and \eqref{eq:91}, we can obtain
\begin{align}
 \frac{1}{T}\sum_{t=1}^T\mathbb{E}\|\mathcal{G}_{\mathcal{X}}(x_t,\nabla F(x_t),\gamma)\|\leq \frac{1}{T} \sum_{t=1}^T \mathbb{E}\big[\mathcal{M}_t\big]  \leq \frac{2\sqrt{3M}m^{1/6}}{T^{1/2}} + \frac{2\sqrt{3M}}{T^{1/3}}.
\end{align}

\end{proof}

\begin{theorem} \label{th:A4}
(Restatement of Theorem 4)
Suppose the sequence $\{x_t,y_t\}_{t=1}^T$ be generated from Algorithm \ref{alg:2}. When $\mathcal{X}=\mathbb{R}^{d_1}$, and  given $B_t=b_tI_{d_2} \ (\hat{b}\geq b_t\geq b >0)$,
 $\eta_t=\frac{k}{(m+t)^{1/3}}$, $\alpha_{t+1}=c_1\eta^2_t$, $\beta_{t+1}=c_2\eta^2_t$, $c_1 \geq \frac{2}{3k^3} + \frac{9\mu^2}{4}$ and $c_2 \geq \frac{2}{3k^3} + \frac{75L^2_f}{2}$, $m \geq \max\big( k^3, (c_1k)^3, (c_2k)^3\big)$, $0< \lambda \leq \min\big(\frac{27\mu bq}{32},\frac{b}{6L_f}\big)$ and
$0< \gamma \leq \min\big( \frac{\rho\lambda\mu\sqrt{q}}{L_f\sqrt{32\lambda^2+150q\kappa^2\hat{b}^2}}, \frac{m^{1/3}\rho}{2Lk} \big)$,
we have
\begin{align}
 \frac{1}{T} \sum_{t=1}^T \mathbb{E}\|\nabla F(x_t)\| \leq \frac{\sqrt{\frac{1}{T}\sum_{t=1}^T\mathbb{E}\|A_t\|^2}}{\rho}\bigg(\frac{2\sqrt{3M'}m^{1/6}}{T^{1/2}} + \frac{2\sqrt{3M'}}{T^{1/3}}\bigg),
\end{align}
where $M' =\frac{\rho(F(x_1) - F^*)}{T\gamma k} + \frac{9L^2_f b_1}{k\lambda\mu}\Delta_1^2  + \frac{2\sigma^2m^{1/3}}{k^2q\mu^2} + \frac{2k^2(c_1^2+c_2^2)\sigma^2}{q\mu^2}\ln(m+T)$.
\end{theorem}

\begin{proof}
According to the above inequality \eqref{eq:78}, we have
 \begin{align} \label{eq:107}
  \frac{1}{T}\sum_{t=1}^T\mathbb{E}\|\nabla F(x_t)\| \leq \frac{1}{T}\sum_{t=1}^T\mathbb{E}\big[\mathcal{M}_t \|A_t\|\big] \leq \sqrt{\frac{1}{T}\sum_{t=1}^T\mathbb{E}[\mathcal{M}_t^2]} \sqrt{\frac{1}{T}\sum_{t=1}^T\mathbb{E}\|A_t\|^2}.
 \end{align}
By using the above inequality \eqref{eq:102} and $ \mathcal{M}_t = \frac{1}{\gamma}\|x_t - \tilde{x}_{t+1}\|
+ \frac{1}{\rho}\big( L_f \|y^*(x_t) - y_t\| + \|\nabla_x f(x_t,y_t)- v_t\|\big)$, we have
\begin{align} \label{eq:108}
 \frac{1}{T}\sum_{t=1}^T\mathbb{E}\big[\mathcal{M}^2_t\big] & \leq \frac{1}{T} \sum_{t=1}^T \mathbb{E}\big[ \frac{3L_f^2}{\rho^2}\|y^*(x_t)-y_t\|^2 + \frac{3}{\rho^2}\|\nabla_{x} f(x_t,y_t) -v_t\|^2 +\frac{3}{\gamma^2}\|\tilde{x}_{t+1}-x_t\|^2 \big] \nonumber \\
 & \leq \frac{12M}{T}(m+T)^{1/3}.
\end{align}
According to the above inequalities \eqref{eq:107} and \eqref{eq:108}, we have
\begin{align}
 \frac{1}{T}\sum_{t=1}^T\mathbb{E}\|\nabla F(x_t)\| \leq \sqrt{\frac{1}{T}\sum_{t=1}^T\mathbb{E}\|A_t\|^2}\frac{2\sqrt{3M}}{T^{1/2}}(m+T)^{1/6}.
\end{align}
Let $M' = \rho^2 M =\frac{\rho(F(x_1) - F^*)}{T\gamma k} + \frac{9L^2_f b_1}{k\lambda\mu}\Delta_1^2  + \frac{2\sigma^2m^{1/3}}{k^2q\mu^2} + \frac{2k^2(c_1^2+c_2^2)\sigma^2}{q\mu^2}\ln(m+T)$,
we have
\begin{align}
\frac{1}{T} \sum_{t=1}^T \mathbb{E}\|\nabla F(x_t)\| \leq \frac{\sqrt{\frac{1}{T}\sum_{t=1}^T\mathbb{E}\|A_t\|^2}}{\rho}\bigg(\frac{2\sqrt{3M'}m^{1/6}}{T^{1/2}} + \frac{2\sqrt{3M'}}{T^{1/3}}\bigg).
\end{align}

\end{proof}

\begin{corollary} (Restatement of Corollary 1)
Under the same conditions of Theorems \ref{th:3} and \ref{th:4}, given mini-batch size $q=O(\kappa^{\nu})$ for $\nu>0$ and $\frac{27\mu bq}{32} \leq \frac{b}{6L_f}$, i.e., $q=\kappa^{\nu} \leq \frac{16}{81L_f\mu}$, our VR-AdaGDA algorithm has a lower gradient complexity of
$\tilde{O}\big(\kappa^{(4.5-\frac{\nu}{2})}\epsilon^{-3}\big)$ for finding an $\epsilon$-stationary point.
\end{corollary}

\begin{proof}
Under the same conditions of Theorems \ref{th:3} and \ref{th:4}, without loss of generality, let $k=O(1)$, $b=O(1)$, $\hat{b}=O(1)$ and $\frac{\rho\lambda\mu\sqrt{q}}{L_f\sqrt{32\lambda^2+150q\kappa^2\hat{b}^2}}\leq \frac{m^{1/3}\rho}{2Lk}$,
we have $m\geq \big(k^3, (c_1k)^3, (c_2k)^3, \frac{8(Lk\lambda\mu)^3q^{3/2}}{L_f(32\lambda^2+150q\kappa^2\hat{b}^2)^{3/2}} \big)$. Let
$\gamma = \frac{\rho\lambda\mu\sqrt{q}}{L_f\sqrt{32\lambda^2+150q\kappa^2\hat{b}^2}}=\frac{\rho\lambda\sqrt{q}}{\kappa\sqrt{32\lambda^2+150q\kappa^2\hat{b}^2}}$
and $\lambda = \min\big(\frac{27\mu bq}{32},\frac{b}{6L_f}\big)$.

Given $q=O(\kappa^{\nu})$ for $\nu>0$ and $\frac{27\mu bq}{32} \leq \frac{b}{6L_f}$, i.e., $\kappa^{\nu}\leq \frac{16}{81L_f\mu}$, it is easily verified that $\lambda=O(q\mu)$, $\gamma=O(\frac{q}{\kappa^3})$, $c_1=O(1)$ and
$c_2=O(L_f^2)$. Due to $L=L_f(1+\kappa)$ and $q \leq \frac{16}{81L_f\mu}$, we have $m=O(L_f^6)$. Then we have
$M=O(\frac{\kappa^3}{q} + \frac{\kappa^2}{q} + \frac{\kappa^2}{q}\ln(m+T)) = O(\frac{\kappa^3}{q})=O(\kappa^{(3-\nu)})$.  Thus, our VR-AdaGDA algorithm has a convergence rate of $O\big(\frac{\kappa^{(3/2-\frac{\nu}{2})}}{T^{1/3}}\big)$.
Let $\frac{\kappa^{(3/2-\frac{\nu}{2})}}{T^{1/3}} \leq \epsilon$, i.e., $\mathbb{E}\big[\mathcal{M}_\zeta\big] \leq \epsilon$ or $\mathbb{E}\|\nabla  F(x_\zeta)\|\leq \epsilon$, we
choose $T \geq  \kappa^{(9/2-\frac{3\nu}{2})}\epsilon^{-3}$.
Thus, our VR-AdaGDA algorithm reaches a lower gradient complexity of $4q \cdot T=O\big(\kappa_y^{(4.5-\frac{\nu}{2})}\epsilon^{-3}\big)$ for finding an $\epsilon$-stationary point.
\end{proof}

\vfill

\end{document}